\documentclass[10pt, reqno]{amsart}

\usepackage{hyperref} 
\usepackage{amssymb, mathrsfs,bbold}
\usepackage{amsmath, amsthm, mathtools}
\usepackage{enumerate}
\usepackage{dsfont}
\usepackage{color}
\usepackage[a4paper]{geometry}
\usepackage[utf8]{inputenc}   
\usepackage{stmaryrd}
\usepackage{titlesec, wasysym}
\usepackage{appendix}
\usepackage{todonotes, fancyhdr}
\usepackage{chngcntr}
\usepackage{cancel}
\usepackage{longtable}

\usepackage{tikz-cd}
\usetikzlibrary{calc}

\usepackage{setspace}
\renewcommand{\baselinestretch}{0.99}

\usepackage[all]{xy}
\numberwithin{subsection}{section}
\numberwithin{subsubsection}{subsection}
\numberwithin{equation}{section} 

\renewcommand{\labelenumi}{\textsf{(\theenumi)}}

\allowdisplaybreaks[4]

\definecolor{macouleur}{rgb}{0.2,0.5,0.8}

\newenvironment{Dem}[1][\unskip]{%
    \begin{list}{\hspace{1.15cm}{\color{black} {\textbf{Proof #1} --}}}{   
        \setlength{\topsep}{0pt}%
        \setlength{\leftmargin}{0pt}%
        \setlength{\rightmargin}{0pt}%
        \setlength{\listparindent}{0pt}%
        \setlength{\itemindent}{0pt}%
        \setlength{\parsep}{0pt}%
        \addtolength{\leftmargin}{0pt} 
        \addtolength{\rightmargin}{0pt}%
    } \item }{\hfill {\color{black} $\rhd$}\end{list}\smallskip}

\newenvironment{Dem*}[1][\unskip]{%
    \begin{list}{\hspace{0cm}{\sf \textbf{{\small Proof #1 --}}}}{   %
        \setlength{\topsep}{0pt}%
        \setlength{\leftmargin}{0pt}%
        \setlength{\rightmargin}{0pt}%
        \setlength{\listparindent}{0pt}%
        \setlength{\itemindent}{0pt}%
        \setlength{\parsep}{0pt}%
        \addtolength{\leftmargin}{20pt}%
        \addtolength{\rightmargin}{0pt}%
    } \item }{\hfill $\rhd$\end{list}\smallskip}

\renewcommand\thesection       {\arabic{section}}
\renewcommand\thesubsection    {\thesection{\boldmath $.$}\arabic{subsection}}
\renewcommand\thesubsubsection    {\thesection{\boldmath $.$}\arabic{subsection}{\boldmath $.$}\arabic{subsubsection}} 

\titleformat{\section}[block] 
{\filcenter\normalfont\sffamily\bfseries\Large}  
{{\hspace{-0.87cm}}{\color{black} \thesection} \hspace{0em} {\color{black} --}\vspace{0cm}}{0.5em}{\color{black} {}} 

\titleformat{\subsection}[runin]
{\filcenter\normalfont\sffamily\bfseries\large}   
{{\hspace{0cm}}{\color{black} \thesubsection} \hspace{0em} {\color{black} --} \vspace{0.1cm}}{.2em}{\color{black} {}}   
\titlespacing{\subsection}{-0pc}{1.5ex plus .1ex minus .2ex}{0pc}   

\titleformat{\subsubsection}[runin]
{\filcenter\normalfont\sffamily\bfseries}   
{\filright\sffamily{\hspace{0cm}}{\color{black} \thesubsubsection}\hspace{0em} {\color{black} --}}{0.4em}{\color{black} {}}\titlespacing{\subsection}{-0pc}{1.5ex plus .1ex minus .2ex}{0pc}

\newtheoremstyle{mystyle}
{3pt}               
{3pt}               
{\it }                      
{}                      
{\bfseries}      
{}                      
{0.5em}                 
{\hspace{0cm}{\color{black} \textrm{#2 --} {\hspace{-0.25cm}} \textrm{ #1}}}   

\theoremstyle{mystyle}

\newtheorem{thm}{Theorem.}   
\newtheorem*{thm*}{Theorem}

\newtheorem{cor}[thm]{Corollary.}
\newtheorem{lem}[thm]{Lemma.}
\newtheorem{prop}[thm]{Proposition.}
\newtheorem*{rem*}{Remark.}

\newtheoremstyle{mystyle3}
{3pt}               
{3pt}               
{\it }                      
{}                      
{\bfseries}      
{}                      
{0.5em}                 
{\hspace{-0.8cm}{\textbf{\textit{#2}} --} {\hspace{-0.02cm}}{\textbf{\textit{#1}}}}

\theoremstyle{mystyle3}

\newtheoremstyle{mystyle2}
{3pt}               
{3pt}               
{\it }                      
{}                      
{\sffamily}    
{}                      
{0.5em}                 
{\llap{#2 }{\it #1{\hspace{0.2cm}--}}}

\theoremstyle{mystyle}

\newtheorem*{definition*}{Definition}
\newtheorem*{theorem*}{Theorem}
\newtheorem*{Remark*}{Remark}
\newtheorem*{lem*} {Lemma}
\newtheorem*{defn*} {Definition}
\newtheorem*{prop*} {Proposition}
\newtheorem*{cor*} {Corollary}

\newcommand{\ssk}{\smallskip}

\renewcommand{\epsilon}{\varepsilon}
\newcommand{\eps}{\epsilon}

\newcommand\bbE{\mathbb{E}}

\newcommand\bbH{\mathbb{H}}
	\newcommand{\ubbH}{\underline{\mathbb{H}}}

\newcommand\bbK{\mathbb{K}}
\newcommand\bbN{\textbf{\textsf{N}}}
\newcommand{\bbP}{\mathbb{P}}
\newcommand\bbQ{\textbf{\textsf{Q}}}
\newcommand\bbR{\textbf{\textsf{R}}}

\newcommand\bbT{\textbf{\textsf{T}}}

\newcommand{\bbOmega}{\mathbb{\Omega}}

\newcommand{\mcE}{\mathcal{E}}
\newcommand{\mcF}{\mathcal{F}}

\newcommand{\mcH}{\mathcal{H}}
\newcommand{\mcI}{\mathcal{I}}

\newcommand{\mcK}{\mathcal{K}}
\newcommand{\mcL}{\mathcal{L}}

\newcommand{\mcQ}{\mathcal{Q}}

\newcommand\mcT{\mathcal T}

\newcommand\mcW{\mathcal W}

\newcommand{\bfB}{\mathbf{B}}

\newcommand{\bfT}{\mathbf{T}}

\newcommand{\bsu}{\boldsymbol{u}}

\newcommand{\sfM}{\mathsf{M}}

\newcommand{\bound}{\textsf{\textbf{bd}}}
\newcommand{\converge}{\textsf{\textbf{cv}}}
\newcommand{\marked}{{\underline{\odot}}}

\newcommand{\scrU}{\mathscr{U}}
\newcommand{\scrV}{\mathscr{V}}
\newcommand{\scrW}{\mathscr{W}}

\newcommand{\res}{\hspace{-0.03cm}:\hspace{-0.03cm}}

\newcommand{\mfs}{\mathfrak{s}}

\newcommand{\id}{\text{\rm id}}

\makeatletter
\newcommand*{\defeq}{\mathrel{\rlap{%
                     \raisebox{0.3ex}{$\m@th\cdot$}}%
                     \raisebox{-0.3ex}{$\m@th\cdot$}}%
                     =}
\makeatother

\makeatletter
\newcommand*{\eqdef}{=\mathrel{\rlap{%
                     \raisebox{0.3ex}{$\m@th\cdot$}}%
                     \raisebox{-0.3ex}{$\m@th\cdot$}}%
                     }
\makeatother

\newcommand{\leby}[1]{%
  \mathrel{\overset{\mathclap{\scriptsize\eqref{#1}}}{\le}}%
}

\newcommand{\spa}{\text{\rm span}}

\begin{document}

\begin{center}
{\LARGE\sffamily{\textbf{Transportation cost inequalities for singular SPDEs}   \vspace{0.5cm}}}
\end{center}

\begin{center}
{\sf I. BAILLEUL, M. HOSHINO} and {\sf R. TAKANO}
\end{center}

\vspace{1cm}

\begin{center}
\begin{minipage}{0.8\textwidth}
\renewcommand\baselinestretch{0.7} \scriptsize \textbf{\textsf{\noindent Abstract.}} We prove that the laws of the BPHZ random models satisfy some transportation cost inequalities in the full subcritical regime if there is no `variance blowup' and the law of the noise is translation invariant and satisfies some transportation cost inequality. We emphasize two consequences of this result or its proof: The automatic integrability properties of the invariant probability measures of a number of singular stochastic partial differential equations, including the $\Phi^4_{4-\delta}$ measures over the $4$-dimensional torus, for all $0<\delta<4$, and a general large deviation principle satisfied by the BPHZ models.   
\end{minipage}
\end{center}

\vspace{0.6cm}

{\it 
\begin{center}
\begin{minipage}[t]{11cm}
\baselineskip =0.35cm
{\scriptsize 

\center{\textbf{\textrm{Contents}}}   \vspace{0.1cm}

	\textrm{1.~Introduction}\dotfill  \pageref{SectionIntro}

	\textrm{2.~Transportation cost inequalities for solutions of singular SPDEs}\dotfill \pageref{SectionTCISPDEs}

	\textrm{3.~Functional setting}\dotfill \pageref{SectionFunctionalSetting}

        \textrm{4.~BPHZ model over a particular regularity-integrability structure}\dotfill \pageref{SectionModelOverregularity-integrability structure}

	\textrm{5.~Transportation cost inequalities for the BPHZ models}\dotfill  \pageref{SectionTCIModels}

	\textrm{6.~Large deviation principle for the BPHZ models}\dotfill  \pageref{SectionLDP}

	\textrm{A.~A relation between transportation cost and spectral gap inequalities}\dotfill  \pageref{SectionAppendix}}
\end{minipage}
\end{center}
}

\vspace{0.8cm}

\section{Introduction}
\label{SectionIntro}

Denote by ${\sf P}(E)$ the set of probability measures on a measurable space $(E,\mcE)$. Fix a measurable function $c:E\times E\rightarrow [0,+\infty]$, usually called a cost function. The $c$-transportation cost between two probability measures $\mu_1,\mu_2$ on $(E,\mcE)$ is defined as
\[
\mcT_c(\mu_1,\mu_2) \defeq \inf \bbE\big[c(X_1,X_2)\big]
\]
for an infimum over the set of all random variables defined on some probability space such that $X_1$ has law $\mu_1$ and $X_2$ has law $\mu_2$. Given $\mu_1,\mu_2$ in ${\sf P}(E)$, the relative entropy $\mcH(\mu_2\vert\mu_1)$ of $\mu_2$ with respect to $\mu_1$ is defined as $\int_E \log\big(\frac{d\mu_2}{d\mu_1}\big) d\mu_2$, if $\mu_2$ is absolutely continuous with respect to $\mu_1$, otherwise we set it equal to $+\infty$. Let $\ell : [0,+\infty]\rightarrow [0,+\infty]$ be a continuous increasing unbounded function with $\ell(0)=0$. A probability measure $\mu$ on $(E,\mcE)$ is said to satisfy an $(\ell,c)$\textit{\textbf{-transportation cost inequality}} if 
\begin{equation} \label{EqTCI}
\ell\big(\mcT_c(\nu,\mu)\big) \leq \mcH(\nu\vert\mu) \qquad (\forall\,\nu\in{\sf P}(E)).
\end{equation}
This type of information theoretic property of a probability measure was first introduced by K. Marton in \cite{Mar86} in relation with the phenomenon of concentration of measure in metric spaces, for $(E,d)$ a metric space, $c(x,y)=d(x,y)^2$ and $\ell(t)$ proportional to $t$. We talk in that case of a $2$-transportation cost inequality. Gozlan \& L\'eonard's review \cite{GL10} provides a nice overview of transportation cost inequalities.

Talagrand \cite{Tal96} proved that Gaussian measures on a Euclidean space satisfy a dimension-free version of \eqref{EqTCI}. This seminal work was followed by numerous other works in different directions. We single out the very influencial work \cite{OV00} of Otto \& Villani which clarified among other things the relation between the $2$-transportation cost inequality and the $\log$-Sobolev and Poincar\'e inequalities, in a Euclidean or Riemannian setting. We also mention Feyel \& \"Ust\"unel's extension of Talagrand's result to an abstract Wiener space setting \cite{FU02}, which opened the door to proving a number of transportation cost inequalities for probability measures that are the laws of some solutions of some stochastic (possibly partial) differential equations, as in the works of Djellout, Guillin \& Wu \cite{DGW04} and Saussereau \cite{Sau12}, or the works of Wu \& Zhang \cite{WZ06}, Khoshnevisan \& Sarantsev \cite{KS19}, Shang \& Zhang \cite{SZ19}, Shang \& Wang \cite{SW20}, Dai \& Li \cite{DL22} or Li \& Wang \cite{LW24}, to cite but a few references.

\ssk

All the preceding results involve some classical It\^o type dynamics. The pathwise solution theories of stochastic dynamics come with some powerful tools that provide some robust and deep informations on the laws of these solutions. The theory of rough paths \cite{Lyons98} gives a pathwise picture of It\^o or Stratonovich stochastic differential equations in which the solution map appears as the composition of a probabilistic lifting map, lifting the noise into a random rough path, and a deterministic solution map, defined on the space of all rough paths and with values in the space of paths. Building on this picture, Riedel proved in \cite{Rie17} that the law of the solution path to a rough differential equation driven by certain classes of Gaussian noises satisfies a transportation cost inequality if the law of the noise satisfies such an inequality. 

\ssk

On the SPDE side, the different pathwise theories developed within the last ten years offer a similar opportunity. It is the purpose of the present work to prove a very general result that provides automatically some transportation cost inequalities for the laws of the solutions of some singular stochastic partial differential equations. 

We will work in the setting of regularity structures \cite{Hai14}. The architecture of that theory is similar to the architecture of rough paths theory. Given any equation in a well identified class of equations, there is a probabilistic lifting map that enriches the noise in the form of a more complex object, called a model. There is also a deterministic, locally Lipschitz, solution map that associates to any model a solution field. Gasteratos \& Jacquier were able to prove in \cite{GJ23} that a transportation cost inequality for some Gaussian noises can be propagated to the solution field of two examples of mildly singular stochastic PDEs, the rough volatility model in finance and the two dimensional parabolic Anderson model equation on the torus. Compared to the rough paths setting, the singular SPDE setting involves an additional difficulty related to the fact that the probabilistic definition of a random model typically involves a non-trivial renormalization procedure. This explains why it is difficult to extend the hand-crafted computations of \cite{GJ23} to some more singular equations, as the renormalization process gets more complex when the equation gets more singular.

\ssk

The map lifting a noise into a model is not unique; we will work here with the Bogoliubov--Parasiuk--Hepp--Zimmermann (BPHZ) prescription first introduced by Bruned, Hairer \& Zambotti in \cite{BHZ}. The lifting map is defined by a limiting procedure whose convergence was first proved in Chandra \& Hairer's unpublished work \cite{ChandraHairer} under a set of assumptions on the noise giving some controls on its cumulants. The limit model is called the BPHZ model. The convergence of the procedure was proved under a different set of assumptions in Hairer \& Steele's work \cite{HS23} and Bailleul \& Hoshino's work \cite{BH23}, where one asks that the law of the noise satisfies a \textit{\textbf{spectral gap}}. Here is what it means. The noise is defined on a Banach space $\Omega$ which contains a dense subspace $H$ which is itself a Hilbert space for some norm $\|\cdot\|_H$ stronger than the norm inherited from $\Omega$. One says that the law of the noise satisfies a spectral gap inequality with constant $a_1>0$ iff one has 
\begin{equation} \label{EqSG}
\bbE\big[(F-\bbE[F])^2\big] \leq a_1 \bbE\big[\Vert dF\Vert_{H^*}^2\big]
\end{equation}
for all cylindrical functions $F(\omega)=f\big(\varphi_1(\omega),\dots,\varphi_n(\omega)\big)$ with $\varphi_i\in \Omega^*$ and some smooth function $f:\bbR^n\rightarrow\bbR$ with at most polynomial growth. The $H$-differential $dF$ of $F$ at point $\omega\in\Omega$ is given for any $h\in H$ by $(d_\omega F)(h)=\sum_{i=1}^n \partial_if\big(\varphi_1(\omega),\dots,\varphi_n(\omega)\big) \varphi_i(h)$. We denoted here by $H^*$ and $\Omega^*$ the topological duals of $H$ and $\Omega$ respectively. We define the cost 
\[
c_H(\omega_1,\omega_2)
=
\begin{cases}
\| h\|_H^2& \text{if}\ \omega_1-\omega_2=h\in H,   \\
+\infty    & \text{otherwise}.
\end{cases}
\]
Different tools were used in \cite{HS23} and \cite{BH23} to prove the convergence of the BPHZ models. Both works follow the inductive strategy introduced first in the seminal work \cite{LOTT} of Linares, Otto, Tempelmayr \& Tsatsoulis to construct the model using a spectral gap assumption. 

\ssk

{\it 1. The main result --} We make here the following stronger assumption.

\ssk

\noindent {\textbf{Assumption A}} -- {\it We are given a separable Banach space $\Omega$ of spacetime functions/distributions which contains a dense subspace $H$ which is itself a Hilbert space for some norm $\|\cdot\|_H$ stronger than the norm inherited from $\Omega$. We equip $\Omega$ with its Borel $\sigma$-algebra. The noise is defined as the identity map on $\Omega$. Its law $\bbP\in {\sf P}(\Omega)$ is centered, invariant under spacetime shifts, and satisfies the transportation cost inequality
\begin{equation} \label{EqTCINoise}
a_\circ \mcT_{c_H}(\mathbb{Q}, \bbP) \leq \mcH(\mathbb{Q} \vert \bbP) \qquad (\forall\,\mathbb{Q} \in {\sf P}(\Omega)).
\end{equation}
for some positive finite constant $a_\circ$.
}

\ssk

In our situation, the space $\Omega$ will be a certain Banach space of distributions on $\bbR^d$. Feyel \& \"Ust\"unel proved in \cite{FU02} that the Gaussian measures on Banach spaces satisfy this assumption with $H$ their Cameron--Martin space. It is a folklore result that the spectral gap inequality \eqref{EqSG} holds under Assumption A, with $a_1=1/(2a_\circ)$. A proof of this result is given in Appendix \ref{SectionAppendix} as we could not find an easily accessible direct proof in the literature. The BPHZ model is thus well defined under Assumption A, using the constructions of \cite{HS23} or \cite{BH23}.

We use in the present work a refinement of the construction of the BPHZ model given in Bailleul \& Hoshino's work \cite{BH23}. The introduction of this refinement is motivated by Gasteratos \& Jacquier's extended contraction principle for transportation cost inequalities; it takes in our setting the following form -- see Lemma 2.11 in \cite{GJ23} for the original statement and Section \ref{SubsectionTCIFacts} for a proof. Let $\textbf{\textsf{M}}$ be a metric space.

\ssk

\begin{prop} \label{PropGeneralisedTCI}
Let $ M: \Omega\rightarrow \textbf{\textsf{M}}$ be a random variable and $c_{\textbf{\textsf{M}}} : \textbf{\textsf{M}}\times\textbf{\textsf{M}}\rightarrow [0,+\infty]$ be a measurable function such that 
\begin{equation} \label{EqConditionExtendedContractionPrinciple}
c_{\textbf{\textsf{M}}}\big(M(\omega_1),M(\omega_2)\big) \leq \mathbb{L}(\omega_2)\max\Big(\sqrt{c_H(\omega_1,\omega_2)} , \sqrt{c_H(\omega_1,\omega_2)}^{1/r}\Big),
\end{equation}
for a finite exponent $r\geq 1$ and all $\omega_1,\omega_2$ in some measurable subset of $\Omega$ of $\bbP$-probability $1$, for $\mathbb{L}\in \bigcap_{2\leq p<\infty} L^p(\Omega,\bbP)$. Assume further that Assumption A holds. Then the law $P_M=\bbP\circ M^{-1}$ of $M$ satisfies for any $1\leq \alpha<2$ the transportation cost inequality 
\begin{equation} \label{EqGasteratosJacquier}
a_\alpha \mcT_{c^\alpha_{\textbf{\textsf{M}}}}(Q , P_M) \leq \mcH(Q\vert P_M)^{\alpha/(2r)} + \mcH(Q\vert P_M)^{\alpha/2}   \qquad   (\forall\,Q\in{\sf P}(\textbf{\textsf{M}}))
\end{equation}
where $a_\alpha = (2\Vert \mathbb{L}^\alpha\Vert_{L^{2/(2-\alpha)}(\bbP)})^{-1}\min(a_\circ^{\alpha/(2r)} , a_\circ^{\alpha/2}) $.
\end{prop}

\ssk

Equation \eqref{EqGasteratosJacquier} corresponds to an $(\ell_\alpha , c^\alpha_{\textbf{\textsf{M}}})$-transportation cost inequality \eqref{EqTCI} with $\ell_\alpha$ the convex and continuous inverse function of the non-negative increasing continuous concave function $(t\in [0,\infty))\mapsto a_\alpha^{-1} (t^{\alpha/(2r)} + t^{\alpha/2})$. Given the statement of Proposition \ref{PropGeneralisedTCI}, it is not surprising that the crux for proving a transportation cost inequality for the BPHZ random model $\widetilde{\sf M}^{\infty;\epsilon}$ lies in proving a pathwise inequality of the following informal form. Here and in the remainder of this introduction we write $\textbf{\textsf{M}}$ for some particular Polish space of models described in detail in Section \ref{SectionModelOverregularity-integrability structure}.

\ssk

\begin{thm} \label{ThmMain}
Under a spectral gap assumption on the law of $\xi$, and a mild condition on the regularity structure considered, one has
\begin{equation} \label{EqTranslationEstimate}
\big\Vert \widetilde{\sf M}^{\infty;\epsilon}(\omega+h) : \widetilde{\sf M}^{\infty;\epsilon}(\omega) \big\Vert_{\textbf{\textsf{M}}} \leq \mathbb{L}(\omega) \max\big(\Vert h\Vert_H , \Vert h\Vert_H^{1/r} \big)
\end{equation}
for all $\omega\in \Omega$ and all $h\in H$, for some exponent $r\geq 1$ and some random variable $\mathbb{L}$ in all the $L^p(\Omega,\bbP)$ spaces $(1\leq p<\infty)$. 
\end{thm}   

\ssk

This result is our main contribution. It will be proved in Theorem \ref{ThmHolderEstimateBPHZ} below, with all the notations defined in the meantime. Proposition \ref{PropGeneralisedTCI} will then entail that the law of the BPHZ model satisfies a family of transportation cost inequalities, from which it will classically follow that its norm has a Gaussian tail. 

\ssk

\begin{cor} \label{CorMainThm}
Under {\bf Assumption A} one has $\bbE\big[\exp(a \Vert \widetilde{\sf M}^{\infty;\epsilon}\Vert_{\textbf{\textsf{M}}}^2)\big] < \infty$ for some positive constant $a$.
\end{cor}   

\ssk

We note here that the norm $ \Vert\cdot\Vert_{ \textbf{\textsf{M}}}$ involved in \eqref{EqTranslationEstimate} and the above expectation is homogeneous with respect to some natural dilation operation on the space of models.

\ssk

We prove Theorem \ref{ThmMain} by constructing the BPHZ model on a particular regularity structure. Like in \cite{LOTT, HS23, BH23} the construction is iterative and uses a spectral gap inequality as a key ingredient to implement the strategy. When formulated with the usual inhomogeneous norm on models, the upper bound \eqref{EqTranslationEstimate} becomes polynomial in $\Vert h\Vert_H$. This fact emphasizes a fundamental difference between the present situation and the constructions in \cite{LOTT, HS23, BH23}. While it was sufficient in these works to use some first order Malliavin calculus, it seems here necessary to use some higher order calculus to get these polynomial estimates. To run this analysis smoothly requires that we work with a regularity structure different from the regularity structures of \cite{HS23} or \cite{BH23}.   

\medskip

{\it 2. Two applications: integrability properties of invariant probability measures and a large deviation principle for the BPHZ model.} One can draw from the statement of Theorem \ref{ThmMain}, or its proof, and Corollary \ref{CorMainThm} a number of consequences. We present two of them here.

\ssk

{\it 2.1 Invariant probability measures --} The question of the longtime existence of some solutions to some singular SPDEs is subtler than the corresponding question for rough differential equations. It is only very recently that some first longtime existence results were obtained for some classes of equations, by Chandra, Feltes \& Weber \cite{CFW24}, Shen, Zhu \& Zhu \cite{SZZ} and Chevyrev \& Gubinelli \cite{ChevyrevGubinelli}. Independently of these examples, under the assumption that a given subcritical singular stochastic PDE is defined on a fixed time interval, it is possible to infer from the transportation cost inequality satisfied by the BPHZ model that the law of the solution field satisfies as well a transportation cost inequality, assuming some a priori pathwise bounds on the solution field. Such a statement will be the object of Corollary \ref{CorTCISolutions} below. Incidentally, Theorem \ref{ThmHolderEstimateBPHZ} and Corollary \ref{CorTCISolutions} allow to extend Riedel's transportation cost result \cite{Rie17} on the law of solutions to some rough differential equations driven by some Gaussian $p$-rough paths with $2<p<3$ to some arbitrary geometric, or branched, $p$-rough paths ($1\leq p<\infty$) whose laws satisfy some transportation cost inequalities. This includes some non-Gaussian examples. 

We obtain as a consequence of Corollary \ref{CorTCISolutions} some transportation cost inequalities for some $\Phi^4$ measures in the full subcritical regime, and more generally for a whole class of invariant probability measures of some Markovian dynamics generated by some singular stochastic PDEs. Let us work over the $4$-dimensional torus $\bbT^4$. For $\nu>\kappa>0$, with $\nu$ irrational and $u_0\in C^{-1-\kappa+\nu}(\bbT^4)$, consider the subcritical singular stochastic PDE
\begin{equation} \label{EqPhi44Minus}
(\partial_t-\Delta)u = - u^3 + (-\Delta)^{-\nu/2}(\zeta)
\end{equation}
where $\zeta$ stands for a spacetime white noise over $\bbT^4$. Denote by $\widetilde{\sf M}^{\infty;\epsilon}$ the BPHZ model over the regularity structure associated with this equation; it satisfies the H\"older bound \eqref{EqTranslationEstimate} for some exponent $r(\nu)$ depending on $\nu$. Equation \eqref{EqPhi44Minus} has a unique solution $u_t(u_0)$ in the setting of regularity structures. It defines a Markovian dynamics on $C^{-1-\kappa+\nu}(\bbT^4)$, defined globally in time and which satisfies the $u_0$-uniform polynomial bound
\[
\Vert u_1(u_0)\Vert_{C^{-1-\kappa+\nu}} \lesssim \big(1 +  \Vert \widetilde{\sf M}^{\infty;\epsilon}\Vert_{\textbf{\textsf{M}}}\big)^{p(\nu)}
\]
for any fixed finite exponent $p(\nu)>1/\nu$. This a priori estimates was proved by Chandra, Moinat \& Weber in \cite{CMW23}. See also Corollary 6.38 of Chevyrev \& Gubinelli \cite{ChevyrevGubinelli}. Denote by $\pi_t(u_0)\in{\sf P}(C^{-1-\kappa+\nu}(\bbT^4))$ the law of $u_t(u_0)$. The bound above entails that the family of probability measures 
\[
\Big(\frac{1}{s}\int_1^{s+1} \pi_t(u_0)dt\Big)_{s\geq 1}
\] 
on $C^{-1-\kappa+\nu}(\bbT^4)$ is tight. Any limit point of this family as $s$ goes to $+\infty$ is an invariant probability measure for the Markovian dynamics \eqref{EqPhi44Minus}. It can be inferred from the deep results of Hairer \& Mattingly \cite{HairerMattingly} on the strong Feller property, and the results of Hairer \& Sch\"onbauer \cite{HairerSchonbauer} on the support theorem for the BPHZ model, that the Markovian dynamics \eqref{EqPhi44Minus} has a unique invariant probability measure. We call such a measure the $\Phi^4_{4-\nu}$ measure and denote it by $\mu$.

Set $\underline{\ell}(t) \defeq t^2$ and define for any $v_1,v_2$ in $C^{-1-\kappa+\nu}(\bbT^4)$ the cost function 
\[
c_\nu(v_1,v_2) = \Vert v_2-v_1\Vert^{1/p(\nu)}_{C^{-1-\kappa+\nu}}.
\]

\ssk

\begin{thm} \label{ThmPhi4Measure}
There is a constant $a\in (0,+\infty)$ such that the $\Phi^4_{4-\nu}$ probability measure on $C^{-1-\kappa+\nu}(\bbT^4)$ satisfies an $(a \underline{\ell},c_\nu)$-transportation cost inequality. Equivalently, there exists a small positive constant $a$ such that
$$
\int_{C^{-1-\kappa+\nu}(\bbT^4)}\exp(a\|\phi\|_{C^{-1-\kappa+\nu}}^{2/p(\nu)})\mu(d\phi) < \infty
$$   
\end{thm}

\ssk

See Theorem \ref{thm:expint} and the proof of Theorem \ref{ThmInvariantMeasures} below. A similar estimate also holds for the $\Phi^4_3$ measure, which somehow corresponds to the case $\nu=1/2$, and then we have the integrability of $\exp(a\|\phi\|^{1-})$ under the $\Phi^4_3$ measure. This fact coincides with the results obtained in Moinat \& Weber \cite{MW20} or Gubinelli \& Hofmanov\'a \cite{GH21}, although they are stated in a somewhat different form. On the other hand, Hairer \& Steele \cite{HS22} obtained a stronger  integrability with a quartic term in the exponent. The relevance of such exponential integrability properties stems from their relation to the regularity axiom in the Osterwalder-Schrader axioms in Euclidean quantum field theory. Our result establishes exponential integrability under the sole assumption that the law of the noise satisfies a certain transportation cost inequality and the solution of the SPDE admits some initial value-free polynomial a priori bound in terms of the model.   

\ssk

Theorem \ref{ThmInvariantMeasures} below is more general than Theorem \ref{ThmPhi4Measure} and shows that a whole class of invariant probability measures of some Markovian dynamics satisfies some explicit transportation cost inequalities. It is expected that the $\Phi^4_{4-\nu}$ measures satisfy some stronger $\log$-Sobolev inequalities. The methods of the present work rely on some pathwise local H\"older estimates of the form \eqref{EqTranslationEstimate}; they do not allow a direct transfer of some $\log$-Sobolev inequality unlike some global Lipschitz estimates. On the other hand, they are robust enough to apply in some situations where it is not clear that the methods used to prove some $\log$-Sobolev inequalities have a chance to be effective, as for instance the methods developed in Bauerschmidt, Bodineau and Dagallier's works \cite{BauerschmidtBodineau,BauerschmidtDagallier}.

\medskip

{\it 2.2 Large deviation principle for the BPHZ model --} Theorem \ref{ThmMain} holds under a spectral gap assumption on the law of the noise. Corollary \ref{CorMainThm}, or rather Corollary \ref{CorTCIModel} below, holds under the stronger assumption {\bf A}. These assumptions go strictly beyond the Gaussian setting as they are robust with respect to perturbations of the measure by a bounded density (Holley-Stroock) or by maps from $\Omega$ into itself that preserve the Hilbert space $H$ of Assumption {\bf A} and are globally Lipschitz on $H$.

Assume further here that $H$ is compactly embedded in the Banach space $\Omega$, that $\xi(\omega)=\omega$ and that the family of random variables $(\kappa\xi)_{0<\kappa\leq 1}$ satisfies a large deviation principle with rate $\kappa^2$ and rate function $I : \Omega\rightarrow[0,+\infty]$ equal to $\|\omega\|_H^2/2$ if $\omega\in H$ and $+\infty$ otherwise. We state here informally the fact that the law of the BPHZ model $\widetilde{\sf M}_\kappa^{\infty;\eps}$ associated with $\kappa\xi$ satisfies then an explicit large deviation principle; the precise statement is given in Theorem \ref{thm:LDPforBPHZmodel} in Section \ref{SectionLDP}.

\ssk

\begin{thm} \label{ThmLDPIntro}
The BPHZ random model $\widetilde{\sf M}_\kappa^{\infty;\eps}$ satisfies a large deviation principle with rate $\kappa^2$ and rate function
$$
J({\sf M}) = \inf\big\{ I(h)\,;\,h\in H,\ L(h)={\sf M} \big\},
$$
where $L(h)\defeq{\sf M}^{h;\eps}$ is the naive model defined by setting ${\sf\Pi}^h(\ocircle)=h$. 
\end{thm}

\ssk

This result extends the seminal large deviation result of Hairer \& Weber \cite{HW15} who proved that result for the particular case of the BPHZ model associated with the Allen-Cahn equation driven by a $(1+1)$-dimensional spacetime white noise, using some Wiener chaos decompositions and some estimate one some Feynman graphs.

\bigskip

\noindent \textit{\textbf{Organisation of the work --}} We dedicate Section \ref{SectionTCISPDEs} to drawing some of the consequences that the fundamental H\"older regularity estimate \eqref{EqTranslationEstimate} has for the laws of the solutions of some singular stochastic PDEs that are well defined globally in time. We prove in particular in Corollary \ref{cor:expintmodel} that the homogeneous norms of the BPHZ models have some Gaussian tails. The above mentioned Theorem \ref{ThmInvariantMeasures} and Theorem \ref{ThmPhi4Measure} are proved in this section. We describe our functional setting in Section \ref{SectionFunctionalSetting}. Section \ref{SectionModelOverregularity-integrability structure} is dedicated to the description of a particular regularity-integrability structure and to the construction of an extended BPHZ model over this setting. This construction is a variation of the construction of the BPHZ model given in \cite{BH23}. We prove that the extended BPHZ model satisfies the H\"older regularity estimate \eqref{EqTranslationEstimate} in Section \ref{SectionTCIModels}. We prove Theorem \ref{ThmLDPIntro} on large deviations in Section \ref{SectionLDP}, as the proof of this result involves a number of technical points involved in the preceding sections. We prove in Appendix \ref{SectionAppendix} that the spectral gap inequality \eqref{EqSG} holds under Assumption A, with $a_1=1/(2a_\circ)$.

\medskip

We refer the reader to Bailleul \& Hoshino's Tourist guide \cite{BHRSGuide} for some background on regularity structures. 

\medskip

\noindent \textbf{\textit{Notations --}} We collect here a few notations that are used throughout the text. For any multiindex $k=(k^j)_{j=1}^{1+d}\in\bbN^{1+d}$ we define
$$
|k|_{\mfs} \defeq 2k_1 + \sum_{j=2}^{1+d} k^j,  \qquad   k! \defeq \prod_{j=1}^{1+d} k^j !,  \qquad  \partial^k\defeq\prod_{j=1}^{1+d}\partial_j^{k^j},
$$
where $\partial_j$ is the partial derivative with respect to the $j$-th variable $x_j$ on $\bbR^{1+d}$. For $k,l\in\bbN^{1+d}$, we write $k\le l$ if $k^j\le l^j$ for any $j\in\{1,\dots,1+d\}$, and we sets $\binom{l}{k} \defeq \prod_{j=1}^{1+d} \binom{l^j}{k^j}$. For every $x=(x_j)_{j=1}^{1+d}\in\bbR^{1+d}$ and $k=(k^j)_{j=1}^{1+d}\in\bbN^{1+d}$ we define
$$
x^k \defeq \prod_{j=1}^{1+d} x_j^{k^j},  \qquad  \|x\|_\mfs \defeq \sqrt{\vert x_1\vert} + \sum_{j=2}^{1+d} |x_j|.
$$


\section{Transportation cost inequalities for solutions of singular SPDEs}
\label{SectionTCISPDEs}

We consider a subcritical singular parabolic stochastic PDEs
\begin{equation} \label{eq:generalSPDE}
(\partial_t-\Delta)u = \mcF(u,\nabla u;\xi)  \qquad  ( t>0,\ x\in\bbT^d)
\end{equation}
with initial condition $u_0\in C^\gamma(\bbT^d)$ for some exponent $\gamma\in\bbR$ and $\xi$ a random (spacetime) noise. The function $\mcF$ is affine in its $\xi$ argument. We will show below that there is a certain extension of the Bruned--Hairer--Zambotti (BHZ) regularity structure \cite{BHZ} associated with the equation over which one can construct an extended version of the BPHZ random model. We will denote this extended version by $\widetilde{\sf M}^{\infty;\eps}$ in the sequel. The extension of the BHZ regularity structure will be described in Section \ref{SectionParticularRIS} and the random model $\widetilde{\sf M}^{\infty;\eps}$ will be constructed in Section \ref{SubsectionBPHZConvergence}. The metric space
\[
\big(\textbf{\textsf{M}}(\scrV_{\varepsilon})_{w_a} , \|\cdot : \cdot\|_{\textbf{\textsf{M}}_{\textrm{hom}}(\scrV_{\varepsilon})_{w_a}}\big)
\] 
in which the BPHZ model $\widetilde{\sf M}^{\infty;\eps}$ takes its values, and all the notations involved, will be explained in Section \ref{SectionModelOverregularity-integrability structure}; they do not matter presently. In this section we will simply denote this space by $(\textbf{\textsf{M}},\Vert\cdot : \cdot\Vert_{\textbf{\textsf{M}}})$. For ${\sf M} \in \textbf{\textsf{M}}$ we set $\Vert{\sf M}\Vert_{\textbf{\textsf{M}}} \defeq \Vert {\bf 0} : {\sf M}\Vert_{\textbf{\textsf{M}}}$ where $\bf 0$ stands for the zero model over the given regularity structure.

\ssk

We will prove in Theorem \ref{ThmHolderEstimateBPHZ} of Section \ref{SectionTCIModels} that if the law of $\xi$ satisfies Assumption A and the condition \eqref{*BoundedVariance} below holds -- the phenomenon called `variance blow-up' in \cite{Hai24} does not occur, then one has for all $\omega\in\Omega$ and $h\in H$ the inequality 
\begin{equation} \label{EqFundamentalEstimateModels}
\big\|{\widetilde{\sfM}}^{\infty;\eps}(\omega+h) :{\widetilde{\sfM}}^{\infty;\eps}(\omega)\big\|_{\textbf{\textsf{M}}} \leq \mathbb{L}(\omega) \max\big(\|h\|_H , \|h\|_H^{1/m}\big)
\end{equation}
for some random variable $\mathbb{L}\in\bigcap_{1\leq p<\infty} L^p(\Omega,\bbP)$ and some integer $m$ that depends only on $\mcF$ and the law of $\xi$. (We note here that \eqref{EqFundamentalEstimateModels} would take a different form if we were using the usual inhomogeneous norm -- see \eqref{mainbound:inhomversion} below.) This regularity estimate of H\"older type is all we need to draw some strong conclusions on the law of the solution to the equation. We obtain them from two elementary facts that we now recall.

\medskip

\subsection{Two useful facts on transportation cost inequalities. \hspace{0.15cm}}
\label{SubsectionTCIFacts}

To make this work self-contained, we give a proof of Proposition \ref{PropGeneralisedTCI} following Gasteratos \& Jacquier's elementary proof of their extended contraction principle -- Lemma 2.11 in \cite{GJ23}.

\ssk

\begin{Dem}[of Proposition \ref{PropGeneralisedTCI}]
Denote by $P_M=\bbP\circ M^{-1}\in{\sf P}(\textbf{\textsf{M}})$ the law of $M$, and pick $Q\in{\sf P}(\textbf{\textsf{M}})$ with $\mcH(Q \vert P_M) < \infty$. Since $\Omega$ is a Polish space, for any $Q\in{\sf P}(\textbf{\textsf{M}})$ with $\mcH(Q \vert P_M) < \infty$ one has 
\begin{equation} \label{EqEntropyImageMeasure}
\mcH(Q \vert P_M) = \inf \Big\{ \mcH(\mathbb{Q} \vert \bbP) \,;\, \mathbb{Q}\in {\sf P}(\Omega),\, Q=\mathbb{Q}\circ M^{-1}\Big\}.
\end{equation}
Write $\Theta(Q,P_M)\subset {\sf P}(\textbf{\textsf{M}}\times \textbf{\textsf{M}})$ for the set of couplings of $Q\in{\sf P}(\textbf{\textsf{M}})$ and $P_M$, and denote by $E$ the expectation operator associated with any of these couplings. Similarly, write $\Theta(\mathbb{Q},\bbP)\subset {\sf P}(\Omega\times\Omega)$ the set of couplings of $\mathbb{Q}\in{\sf P}(\Omega)$ and $\bbP$ and denote by $\bbE$ the expectation operator associated with any of these couplings. For all $1\leq \alpha<2$, one has $2r/(2r-\alpha)\leq 2/(2-\alpha)$ and 
\begin{equation*} \begin{split}
\inf_{\Theta(Q,P_M)} \, E\big[c_{\textbf{\textsf{M}}}^\alpha(\cdot,\cdot)\big] 
&\leq \inf_{\Theta(\mathbb{Q},\bbP)} \,  \bbE\big[ c_{\textbf{\textsf{M}}}^\alpha(M\cdot , M\cdot) \big]   \\
&\leby{EqConditionExtendedContractionPrinciple}
\inf_{\Theta(\mathbb{Q},\bbP)} \Big( \bbE\big[ \mathbb{L}^\alpha c_H^{\alpha/2}(\cdot , \cdot) {\bf 1}_{c_H\geq 1}\big] + \bbE\big[ \mathbb{L}^\alpha c_H^{\alpha/(2r)}(\cdot , \cdot) {\bf 1}_{c_H<1}\big] \Big)   \\
\end{split} \end{equation*}
\begin{equation*} \begin{split}
&\leq \inf_{\Theta(\mathbb{Q},\bbP)} \Big(\Vert \mathbb{L}^\alpha\Vert_{L^{2/(2-\alpha)}(\bbP)} \bbE[c_H]^{\alpha/2} + \Vert \mathbb{L}^\alpha\Vert_{L^{2r/(2r-\alpha)}(\bbP)} \bbE[c_H]^{\alpha/(2r)} \Big)   \\
&\leby{EqTCINoise}
\Vert \mathbb{L}^\alpha\Vert_{L^{2/(2-\alpha)}(\bbP)} \bigg\{ \Big(\frac{\mcH(\mathbb{Q} \vert \bbP)}{a_\circ}\Big)^{\frac{\alpha}{2}} + \Big(\frac{\mcH(\mathbb{Q} \vert \bbP)}{a_\circ}\Big)^{\frac{\alpha}{2r}} \bigg\}   \\
&\leq b_\alpha \max\big( \mcH(\mathbb{Q} \vert \bbP)^{\alpha/2} , \mcH(\mathbb{Q} \vert \bbP)^{\alpha/(2r)} \big)
\end{split} \end{equation*}
with $b_\alpha = 2 \Vert \mathbb{L}^\alpha\Vert_{L^{2/(2-\alpha)}(\bbP)} \max\big( a_\circ^{-\alpha/(2r)} , a_\circ^{-\alpha/2}\big)$, using H\"older inequality in the third inequality. One then gets the conclusion from \eqref{EqEntropyImageMeasure}.
\end{Dem}

\ssk

The proof makes clear what conclusion can be obtained if we only assume in Proposition \ref{PropGeneralisedTCI} that $\mathbb{L}\in L^p(\Omega,\bbP)$ for some finite $p\geq 2$. The following fact follows from Proposition \ref{PropGeneralisedTCI} and the estimate \eqref{EqFundamentalEstimateModels}. 

\ssk

\begin{cor} \label{CorTCIModel}
The law $\widetilde{P}^{\infty;\epsilon} \in {\sf P}(\textbf{\textsf{M}})$ of the \emph{BPHZ} model $\widetilde{\sf M}^{\infty;\eps}$ satisfies for each $1\leq \alpha<2$ an $(\ell_\alpha , \|\cdot : \cdot\|_{\textbf{\textsf{M}}}^\alpha)$-transportation cost inequality, where $\ell_\alpha$ stands for the convex and continuous inverse function of the non-negative increasing continuous concave function $(t\in [0,\infty))\mapsto a_\alpha^{-1} (t^{\alpha/(2r)} + t^{\alpha/2})$, and the random variable $\mathbb{L}$ involved in the definition of the constant $a_\alpha$ is the one that appears in \eqref{EqFundamentalEstimateModels}.
\end{cor}

\ssk

In the particular case where the cost function is the distance function of a Polish space, one can relate some transportation cost inequalities to some exponential integrability properties of the reference probability measures. The following statement was proved by Gozlan in Theorem $1.13$ of \cite{Goz06}.

\ssk

\begin{thm}[{\cite[Theorem 1.13]{Goz06}}] \label{thm:expint}
Let $\ell : [0,\infty)\to[0,\infty)$ be a convex, continuous increasing function such that $\ell(0)=0$. Its conjugate function $\ell^*:[0,\infty)\to[0,\infty]$ is defined for $0\leq s<\infty$ as
$$
\ell^*(s) = \sup_{r\geq 0}\big\{ sr - \ell(r) \big\}.
$$
Let $(E,d)$ be a Polish space and $\mu$ be a probability measure on $E$. Assume that 
\begin{itemize}
	\item[--] one has $\{\ell^*<\infty\}=[0,s_0)$ for some $s_0>0$, and there are some positive constants $\rho$ and $s_1<s_0$ such that $\ell^*(s)\geq \rho s^2$ for any $s\in[0,s_1]$;
	
	\item[--] there is a reference point $x_0\in E$ such that $\int_E d(x_0,x)\mu(dx)<\infty$.
\end{itemize}	 
In this setting, the following properties are equivalent.
\begin{itemize}
	\item[(a)] The probability $\mu$ satisfies an $\big(\ell(\frac{\cdot}{\kappa_1}),d(\cdot,\cdot)\big)$-transportation cost inequality for some constant $0<\kappa_1<\infty$.   \vspace{0.1cm}
	
	\item[(b)] The integral $\int_E e^{\ell(\kappa_2 d(x_0,x))}\mu(dx)$ is finite for some constant $0<\kappa_2<\infty$.
\end{itemize}
Moreover, if (a) holds true for some $\kappa_1$, one has $\|d(x_0,\cdot)-\int_Ed(x_0,\cdot)d\mu\|_{\ell}\le3\kappa_1$, where $\|\cdot\|_\ell$ is Luxemburg norm defined by
$$
\|f\|_\ell=\inf\Big\{\lambda>0\,;\,\int_E\exp\Big(\ell\Big(\frac{f}\lambda\Big)\Big) d\mu \le 2 \Big\}.
$$
\end{thm}

\ssk

The statement about the Luwemburg norm is not stated explicitly in Gozlan's theorem but it is contained in its proof. We deduce from Corollary \ref{CorTCIModel} with $\alpha=1$ and the (a) $\Rightarrow$ (b) part of Theorem \ref{thm:expint} that $\|\widetilde{\sfM}^{\infty;\varepsilon}\|_{\textbf{\textsf{M}}}$ has a Gaussian tail.

\ssk

\begin{cor} \label{cor:expintmodel}
There is a finite positive constant $\kappa_3$ such that 
\[
\bbE\big[\exp\big(\kappa_3\|\widetilde{\sfM}^{\infty;\varepsilon}\|_{\textbf{\textsf{M}}}^2\big)\big] < \infty.
\] 
\end{cor}

\ssk

\begin{Dem}
We work here with $\alpha=1$ in the conclusion of Corollary \ref{CorTCIModel}. It is straightforward to show that the concave function  $\ell_1$ satisfies the conditions of Theorem \ref{thm:expint}. By applying this theorem in the case where $(E,d)=\big(\textbf{\textsf{M}} , \|\cdot\res\cdot\|_{\textbf{\textsf{M}}}\big)$ and the reference point in $E$ is the zero model, we obtain the existence of a constant $\kappa_2$ such that $\bbE\big[\exp\big(\ell_1\big(\kappa_2\|\widetilde{\sfM}^{\infty;\varepsilon}\|_{\textbf{\textsf{M}}}\big)\big)\big]$ is finite. The result then follows from the fact that $\ell_1(t)\gtrsim t^2$ for $t\ge 1$.

One can also deduce this statement from \eqref{EqFundamentalEstimateModels} and Theorem 3.1 of Riedel's work \cite{Rie17}.
\end{Dem}

\ssk

Any real-valued $1$-Lipschitz function on $\big(\textbf{\textsf{M}} , \|\cdot\res\cdot\|_{\textbf{\textsf{M}}}\big)$ is thus strongly concentrated around its mean
$$
\bbP(f-\bbE[f] \geq t) \leq \exp\hspace{-0.07cm}\big(\hspace{-0.1cm}-ct^2/2\big)
$$
and one has for all $\lambda\in\bbR$
$$
\bbE\big[ \exp(\lambda(f-\bbE[f])) \big] \leq \exp\big(\lambda^2/(2c)\big)
$$
for some positive finite constant $c$. We also deduce the following fact from Corollary \ref{cor:expintmodel} and the (b) $\Rightarrow$ (a) part of Theorem \ref{thm:expint}. Recall that $\underline{\ell}(t)=t^2$.

\ssk

\begin{cor} \label{CorNormControlX}
Let $(E,\Vert\cdot\Vert_E)$ be a Banach space and $X : \Omega\rightarrow E$ be a random variable that satisfies the upper bound
\[
g(\Vert X\Vert_E) \leq 1 + \big\| \widetilde{\sf M}^{\infty;\epsilon} \big\|_{\textbf{\textsf{M}}}
\]
for some continuous concave strictly increasing function $g:[0,\infty)\rightarrow[0,\infty)$ with $g(0)=0$. Then there is a constant $0<a<\infty$ such that the law of $X$ satisfies an $\big(a\underline{\ell} , g(\|\cdot-\cdot\|_E)\big)$-transportation cost inequality.
\end{cor}

\ssk

Note that $g(\|\cdot-\cdot\|_E)$ is a metric on $E$ and $\big(E , g(\|\cdot-\cdot\|_E)\big)$ is a Polish space. In a typical situation $g(\cdot)$ would be the inverse function of a constant multiple of one of the functions $z^p$ or $\exp\big(\kappa z^p\big)-1$, for some constants $1\leq p<\infty$ and $0<\kappa<\infty$.

\medskip

\subsection{Transportation cost inequalities for solutions of singular SPDEs. \hspace{0.1cm}}
\label{SectionTCISolutionSPDEs}

We consider the subcritical singular parabolic stochastic PDEs \eqref{eq:generalSPDE}. We refer to Hairer's original work \cite{Hai14} or Bailleul \& Hoshino's Tourist guide \cite{BHRSGuide} for some background on the regularity structure formulation of Equation \eqref{eq:generalSPDE} in a space of modelled distributions depending on $\widetilde{\sfM}^{\infty;\varepsilon}(\omega)$, and the fact that this reformulation is locally well-posed on a random time interval. This conclusion holds under some regularity assumptions on $\mcF$ that are irrelevant here. We denote by bold $\bsu$ the solution modelled distribution and by $u={\sf R}^{\widetilde{\sf M}^{\infty;\epsilon}}(\bsu)$ its reconstruction via the reconstruction operator associated with the BPHZ model $\widetilde{\sf M}^{\infty;\epsilon}$.

One needs some assumptions on $\mcF$ to ensure the longtime existence of $\bsu$ and $u$. Some systematic investigations for finding such conditions were only done recently in the works of Chandra, Feltes \& Weber \cite{CFW24}, Shen, Zhu \& Zhu \cite{SZZ} and Chevyrev \& Gubinelli \cite{ChevyrevGubinelli}. The first two works deal with some mildly singular equations for which $\mcF(u,\nabla u;\xi)=\sigma(u)\xi$ with $\sigma\in C^2_b(\bbR)$ and $\xi$ has almost surely some regularity comparable to the regularity of the two dimensional space white noise. The third work deals with some $\Phi^4$-like dynamics where $\mcF(u,\nabla u;\xi) = P(u,\nabla u) + f(u,\nabla u)\xi$ and $P$ is coercive in some sense, as in the $\Phi^4$ example where $P(u,\nabla u)=-u^3$ and $f$ is constant. They can also deal with some non-constant functions $f$, depending on the almost sure regularity of $\xi$ -- see Section 6.3 of \cite{ChevyrevGubinelli}. We note that unlike \cite{CFW24,SZZ}, the results of \cite{ChevyrevGubinelli} hold in the full subcritical regime and for some dynamics set in the full Euclidean space.

\ssk

\noindent {\it We assume for the remainder of this section that the singular SPDE \eqref{eq:generalSPDE} is globally well-posed; we denote by $u$ the solution.}

\ssk

In all the situations considered in \cite{CFW24,SZZ,ChevyrevGubinelli} one gets the longtime existence of $\bsu$ and $u$ from some a priori estimates that control the size of $u$ in some spacetime domains in terms of the (BPHZ) model. Corollary \ref{CorNormControlX} then takes the following form, again with $\underline{\ell}(t)=t^2$.

\ssk

\begin{cor} \label{CorTCISolutions}
Assume that the random solution $u$ to Equation \eqref{eq:generalSPDE} takes its values in a Banach space $(E,\Vert\cdot\Vert_E)$ and one has
$$
g(\|u\|_E) \leq 1+\|\widetilde{\sf M}^{\infty;\epsilon}\|_{\textbf{\textsf{M}}}
$$
for some continuous concave strictly increasing function $g:[0,\infty)\rightarrow[0,\infty)$ with $g(0)=0$. Then the law of $u$ satisfies an $\big(a\underline{\ell} , g(\|\cdot-\cdot\|_E)\big)$-transportation cost inequality for some finite positive constant $a$.
\end{cor}

\ssk

The function space $E$ is typically of the form $C\big((0,T];C^\gamma(\bbT^d)\big)$ for some fixed positive finite time $T$, some $\gamma\in \bbR$ and some norm of the form $\sup_{0<t\leq T} t^q \Vert u(t)\Vert_{C^\gamma}$ with $0<q<\infty$. One has for instance $g(z)=z^{1/p}$ for some constant $1\leq p<\infty$ for the solution of the dynamical $\Phi^4_3$ equation \cite{CMW23} and the two-dimensional parabolic Anderson model \cite{CFW24}. To emphasize the dependence of $u$ on the initial condition $u_0$, we denote by $u(t;u_0)$ the value at time $t$ of $u$. We denote below by $(\theta_s)_{s\in[0,\infty)}$ a family of shift operators acting on time-space functions/distributions via the formula $\theta_s(\cdot)(r,x)=(\cdot)(r+s,x)$. In our setting, the probability space $\Omega$ is a distribution space and the shifts act on $\Omega$. Recall from Assumption A that the probability $\bbP$ is invariant by the action of the shifts.

\ssk

\begin{thm} \label{ThmInvariantMeasures}
Assume that there exists an open interval $(\alpha,\beta)\subset\bbR$ such that the dynamics \eqref{eq:generalSPDE} defines a Feller process on $C^\gamma(\bbT^d)$ for any $\gamma\in(\alpha,\beta)$. Assume as well the existence of a continuous concave strictly increasing unbounded function $g:[0,\infty)\rightarrow[0,\infty)$, with $g(0)=0$, such that $u(1;u_0)$ satisfies the $u_0$-uniform bound
\begin{equation} \label{EqComingDownLike}
g(\|u(1;u_0)\|_{C^\gamma}) \leq 1+\|\widetilde{\sfM}^{\infty;\eps}\|_{\textbf{\textsf{M}}}.
\end{equation}
Then the dynamics \eqref{eq:generalSPDE} on $C^\gamma(\bbT^d)$ has an invariant probability measure which further satisfies an $\big(a\underline{\ell},g(\|\cdot-\cdot\|_{C^\gamma})\big)$-transportation cost inequality for some constant $0<a<\infty$.
\end{thm}

\ssk

The condition \eqref{EqComingDownLike} is of course reminiscent of the `coming down from infinity' property satisfied by the solution of the $\Phi^4_3$ dynamics.

\ssk

\begin{Dem}
By the Markov property, we have for any $t\geq 1$
$$
g\big(\|u(t;u_0)\|_{C^\gamma}\big) = g\big(\|u(1;u(t-1))\|_{C^\gamma}\big) \leq 1+\|\theta_{t-1}\widetilde{\sfM}^{\infty;\varepsilon}\|_{\textbf{\textsf{M}}}.
$$
It follows form the shift invariance of $\bbP$ that the law of $\widetilde{\sfM}^{\infty;\varepsilon}$ is also shift-invariant. We have as a consequence
\begin{equation} \label{EqUniformExpBound}
\sup_{1\leq t<\infty} \bbE\big[\exp\big(\kappa_2 g(\|u(t;u_0)\|_{C^\gamma})^2\big) \big] < \infty
\end{equation}
for some positive constant $\kappa_2$. Denote by $\pi_t(u_0)\in{\sf P}\big(C^\gamma(\bbT^d)\big)$ the law of $u(t;u_0)$ and set $\overline{\pi}_s(u_0) \defeq \frac{1}{s} \int_1^{s+1} \pi_t(u_0)dt$ for $s\geq 1$. The uniform bound \eqref{EqUniformExpBound} classically implies the weak convergence of $\overline{\pi}_s(u_0)$ in ${\sf P}\big(C^\gamma(\bbT^d)\big)$ to a probability measure $\mu$ invariant for the dynamics \eqref{eq:generalSPDE}, as $s$ moves along a particular sequence $(s_n)_{n\geq 1}$ diverging to infinity. The reasoning goes as follows:  Define for any constant $0<m<\infty$ the set $K_m = \big\{ f\in C^{\gamma+\delta}(\bbT^d) \,;\, \Vert f\Vert_{C^{\gamma+\delta}} \le m \big\} \subset C^{\gamma}(\bbT^d)$, for $\delta>0$ small. The set $K_m$ is bounded in $C^{\gamma+\delta}(\bbT^d)$ hence compact in $C^{\gamma}(\bbT^d)$. We deduce from the $t$-uniform bound \eqref{EqUniformExpBound} with $\gamma$ replaced by $\gamma+\delta$ implies by Chebychev inequality the upper bound $\overline\pi_s(u_0)(C^{\gamma}(\bbT^d)\setminus K_m) \lesssim e^{-\kappa_2g(m)^2}$, from which the tightness of the $(\overline{\pi}_s(u_0))_{1\leq s<\infty}$ follows since $g$ is increasing with $g(+\infty)=+\infty$. 

The probability measure $\mu$ satisfies
\begin{align*}
\int_{C^\gamma(\bbT^d)} e^{\kappa_2g(\|\phi\|_{C^\gamma})^2} \mu(d\phi)
&\le\liminf_{n\to\infty}\int_{C^\gamma(\bbT^d)} e^{\kappa_2g(\|\phi\|_{C^\gamma})^2} \overline{\pi}_{s_n}(u_0)(d\phi)   \\
&\le\liminf_{n\to\infty}\frac1{s_n}\int_1^{s_n+1}\bigg(\int_{C^\gamma(\bbT^d)} e^{\kappa_2g(\|\phi\|_{C^\gamma})^2} \pi_t(u_0)(d\phi)\bigg)dt   \\
&\leq \sup_{1\leq t<\infty} \bbE\big[\exp\big(\kappa_2 g(\|u(t;u_0)\|_{C^\gamma})^2\big) \big] < \infty.
\end{align*}
The conclusion of the statement follows from the (b) $\Rightarrow$ (a) part of Theorem \ref{thm:expint}.
\end{Dem}

\ssk

As for the $\Phi^4_{4-\nu}$ measure, the uniqueness of an invariant probability measure from the dynamics is a consequence of the fact that the transition semigroup of the Markovian dynamics on $C^\gamma(\bbT^d)$ has the strong Feller property, a fact proved in great generality by Hairer \& Mattingly in \cite{HairerMattingly}, and the support theorem for the BPHZ model proved by Hairer \& Sch\"onbauer in \cite{HairerSchonbauer}. The latter implies that the law of $u$ on the interval $[0,T]$ has support in $C([0,T],C^\gamma(\bbT^d))$ given by all functions with values $u_0$ at time $0$. (This is Theorem 1.12 in \cite{HairerSchonbauer}.) The reasoning in the proof of Theorem 1.13 in \cite{HairerSchonbauer} applies and gives the uniqueness of an invariant measure. Therefore we obtain Theorem \ref{ThmPhi4Measure} on the $\Phi^4_{4-\nu}$ measure as a direct corollary of Theorem \ref{ThmInvariantMeasures} and the polynomial estimates on the solution $u$ to \eqref{EqPhi44Minus} proved by Chandra, Moinat \& Weber in \cite{CMW23}. Consult \cite{ChevyrevGubinelli} for a number of other examples where the reasoning of the proof of Theorem \ref{ThmInvariantMeasures} can be applied and gives some transportation cost inequalities for some invariant probability measures of some Feller dynamics.

\medskip

The next three sections are dedicated to proving our main result, Theorem \ref{ThmMain}. First we introduce in Section \ref{SectionFunctionalSetting} an appropriate functional setting. Our proof of Theorem \ref{ThmMain} requires that we introduce a regularity structure different from the structures used in \cite{HS23} and \cite{BH23}. This is related to the fact that we aim at proving some polynomial estimates which require the use of some high order Malliavin calculus rather than only a first order calculus. This regularity structure is introduced in Section \ref{SectionParticularRIS}. We use the flexibility of the notion of regularity-integrability structure from \cite{BH23} rather than the usual notion of regularity structure. The construction of the BPHZ model over this structure is described in Section \ref{SubsectionBPHZConvergence}. As this is only a variation on the arguments from \cite{BH23} we only describe the structure of the argument and emphasize the differences with the two settings. The proof of the H\"older estimate \eqref{EqTranslationEstimate} itself is given in Section \ref{SectionTCIModels}. 

The mechanics of the proof is simple and can be informally illustrated on the ${\sf \Pi}$-part of the limit model $\widetilde{\sf M}^{\infty,\epsilon}$. This model actually comes under the form of the restriction of a model with parameters in $(\omega,h)\in\Omega\times H$ which we denote here by $\overline{\sf M}^{\omega,h}$. The core of the argument rests on the following type of identities
$$
\overline{\sf \Pi}_x^{\omega+h,0}(\tau) - \overline{\sf \Pi}_x^{\omega,0}(\tau) = \sum_V \overline{\sf \Pi}^{\omega,h}_x(D_V\tau)
$$
where $V$ is some finite index set and $D_V$ is some type of high order derivative operator. The quantities $\overline{\sf \Pi}^{\omega,h}_x(D_V\tau)$ are $|V|$-multilinear in $h$ so they have some estimates proportional to $\Vert h\Vert_H^{|V|}$.

\ssk

\section{Functional setting}
\label{SectionFunctionalSetting}

We describe in this short section the functional setting involved in the sequel -- see Section 2 of \cite{BH23} and references therein for the details. Set 
$$
\mcL \defeq \partial_1^2 - \big(1-\partial_2^2-\cdots-\partial_{1+d}^2\big)^2. 
$$
Its associated semigroup $(\mcQ_s \defeq e^{s\mcL})_{s>0}$ gives a representation of the inverse of the heat operator $\partial_t-\Delta+1$ by
$$
\big(\partial_1-\Delta+1\big)^{-1} = -\int_0^\infty (\partial_1+\Delta-1) \, e^{s\mcL} \, ds.
$$
Denote by $Q_s$ the smooth integral kernel of $\mcQ_s$
$$
(\mcQ_sf)(x) =\mcQ_t(x,f) = \int_{\bbR^{1+d}}Q_s(x-y)f(y)dy.
$$
We define the family of weight functions $(w_a)_{a\ge0}$ on $\bbR^{1+d}$ by
$$
w_a(x) \defeq \big(1+\|x\|_\mfs\big)^{-a}.
$$
For $a\ge0$ and $p\in[1,\infty]$ we define the weighted $L^p$ norm by
$$
\|f\|_{L^p(w_a)}\defeq\|fw_a\|_{L^p(\bbR^{1+d})}.
$$
For any non-negative integer $m$ we set $\mcQ_t^{(m)}=(-t\mcL)^m \mcQ_t$. For every $r<4m$ and $p,q\in[1,\infty]$ we define the \textbf{\textit{Besov space}} $B_{p,q}^{r,\mcQ}(w_a)$ as the completion of $C(\bbR^{1+d})\cap L^p(w_a)$ under the norm
\begin{equation*} 
\|f\|_{B_{p,q}^{r,\mcQ}(w_a)}^{(m)}\defeq\|\mcQ_1f\|_{L^p(w_a)}+\big\|t^{-r/4}\|\mcQ_t^{(m)}f\|_{L^p(w_a)}\big\|_{L^q((0,1];\frac{dt}t)}.
\end{equation*}
The topological space $B_{p,q}^{r,\mcQ}(w_a)$ is defined independently to the choice of $m$ as long as $m>r/4$, since the norms $\|\cdot\|_{B_{p,q}^{r,\mcQ}(w_a)}^{(m_1)}$ and $\|\cdot\|_{B_{p,q}^{r,\mcQ}(w_a)}^{(m_2)}$ are equivalent for $m_1,m_2>r/4$. We set
$$
H^{r,\mcQ}(w_a) \defeq B_{2,2}^{r,\mcQ}(w_a), \qquad C^{r,\mcQ}(w_a) \defeq B_{\infty,\infty}^{r,\mcQ}(w_a).
$$
We define a family $(\widetilde{K}_s)_{0<s\leq 1}$ of $2$-regularizing operators setting 
\[
\widetilde{K}_s = - (\partial_1+\Delta-1)e^{s(\partial_1^2-(\Delta-1)^2)}.
\] 
For technical reasons we fix a compactly supported smooth function $\chi$ which is equal to $1$ on a neighborhood of $0$ and we consider the modified operators
$$
K_t = \big( 1-\chi(\partial) \big) \widetilde{K}_t,
$$
where 
$$
\chi(\partial)f \defeq \mathscr{F}^{-1}\big(\chi \, \mathscr{F} f\big)
$$ 
is a Fourier multiplier operator defined by Fourier transform $\mathscr{F}$ and its inverse $\mathscr{F}^{-1}$. This modification ensures that
$$
\int_{\bbR^d}x^k\partial^lK_t(x)dx = 0
$$
for any $k,l\in\bbN^{1+d}$, where we denote the operator $K_t$ and its integral kernel $K_t$ by the same symbol. 
Finally, for any $f\in C(\bbR^{1+d})\cap L^p(w_a)$ and $k\in\bbN^{1+d}$, we define
\begin{equation} \label{EqQSpacetimeOperator}
\partial^k\mcK(x,f) \defeq \int_0^1 \hspace{-0.15cm} \int_{\bbR^{1+d}} \partial_x^k K_s(x-y)f(y) \, dyds.
\end{equation}
The operator $\mcK$ coincides with $(\partial_1-\Delta+1)^{-1}$ modulo some regularizing operators.

\medskip

\section{BPHZ model over a particular regularity-integrability structure}
\label{SectionModelOverregularity-integrability structure}

We will prove in Theorem \ref{ThmHolderEstimateBPHZ} in Section \ref{SectionTCIModels} the local H\"older regularity estimate \eqref{EqFundamentalEstimateModels} of the BPHZ map.

The BPHZ model was constructed in \cite{HS23} and \cite{BH23} using two different functional settings and different sets of tools. We use here the apparatus of \cite{BH23}, as a variation on the arguments used therein happens to lead directly to \eqref{EqFundamentalEstimateModels}.

A notion of regularity-integrability structure was introduced in \cite{SemigroupMasato} and \cite{BH23} to implement an iterative construction of the BPHZ model involving the Malliavin derivative operator and a spectral gap assumption on the law of the noise. The index set of a regularity structure is a subset of $\bbR$; the index set of a regularity-integrability structure is a subset of $\bbR\times [1,+\infty]$. This is related to the fact that we quantify the regularity of some analytical objects associated with any symbol $\tau$ of the structure in some $\tau$-dependent Besov space $B^{r(\tau)}_{i(\tau),\infty}$ rather than in a fixed function space.

\ssk

We define a strict partial order on $\bbR\times[1,\infty]$ by setting
$$
(r,i)<(s,j) \quad \overset{\text{\rm def}}{\Longleftrightarrow} \quad \Big\{r<s \ \textrm{and} \ \frac{1}{i}\leq \frac{1}{j}\Big\}.
$$
We recall from \cite{SemigroupMasato} and \cite{BH23} that a \textit{\textbf{regularity-integrability structure}} $(A,T,G)$ consists of the following elements.
\begin{enumerate} \setlength{\itemsep}{0.1cm}
\renewcommand{\labelenumi}{\textrm{\bf (\alph{enumi})}}
	\item[--] The index set $A$ is a subset of $\bbR\times[1,\infty]$ such that for every $(s,j)\in\bbR\times[1,\infty]$ the set $\{(r,i)\in A\, ;\, (r,i)<(s,j)\}$ is finite.

	\item[--] The vector space $T=\bigoplus_{{\bf a}\in A}T_{\bf a}$ is an algebraic sum of Banach spaces $(T_{\bf a},\|\cdot\|_{\bf a})$.  

	\item[--] The structure group $G$ is a group of continuous linear operators on $T$ such that one has for all $\Gamma\in G$ and ${\bf a}\in A$
$$
(\Gamma - \textrm{Id})(T_{\bf a})\subset\bigoplus_{{\bf b}\in A,\, {\bf b}<{\bf a}} T_{\bf b}.
$$
\end{enumerate}

\ssk

We introduce in Section \ref{SectionParticularRIS} a particular regularity-integrability structure. This section is dedicated to proving that one can define the BPHZ model over this extended regularity-integrability structure assuming that the law of the noise satisfies Assumption A. We assume from now on that the Hilbert space $H$ involved in that assumption is the Hilbert space $H^{-s_0,\mcQ}(w_b)$, for some $s_0\in\bbR$ and $b\ge0$ which will be specified in Theorem \ref{ThmConstructionBPHZ} below.

\medskip

\subsection{A particular regularity-integrability structure. \hspace{0.15cm}}
\label{SectionParticularRIS}

We introduce in this section a regularity-integrability structure that differs from the regularity-integrability structure of \cite{BH23} by the fact that 
\begin{itemize}
	\item[--] we include some trees corresponding to all order Malliavin derivatives, while only first order derivatives are considered in \cite{BH23},
	\item[--]  we introduce two copies of the $H$ space in our type set
	$$
	{\frak{L}} = \big\{ \bbK, \bbOmega,\bbH,\ubbH \big\}.
	$$
	The type $\bbH$ will represent some elements of $H$ estimated in some $B_{\infty,\infty}$-type Besov norm, while the type $\ubbH$ will represent some elements of $H$ estimated in some $B_{p,\infty}$-type Besov norm, with $p$ varying freely in the interval $[2,\infty]$.
\end{itemize}
We invite the reader familiar with the constructions of the trees from \cite{BHZ} to skip the next paragraph.

\medskip

\noindent \textbf{\S1. Sets of trees --} We denote by $\widetilde{\bf T}$ the set of all 4-tuples $(\tau,{\frak{t}},{\frak{n}},{\frak{e}})$, with $\tau$ a non-planar rooted tree $\tau$ with vertex/node set $N_\tau$ and edge set $E_\tau$, with ${\frak{t}}:E_\tau\to{\frak{L}}$ a label map, and with two decoration maps ${\frak{n}}:N_\tau\to\bbN^{1+d}$ and ${\frak{e}}:E_\tau\to\bbN^{1+d}$. We sometimes abuse notations and write $\tau$ for a decorated tree $(\tau,{\frak{t}},{\frak{n}},{\frak{e}})$ if the context makes clear what the decoration is. A {\it leaf} in a decorated tree is an edge $e=(u,v)$ such that there is no outgoing edge from $v$. (Every edge is represented as an ordered pair $(u,v)$ of two nodes, where $u$ is nearer to the root than $v$.)

Denote by $\widetilde{T}$ the linear space spanned by $\widetilde{\bf T}$. The tree product of $\tau,\sigma\in\widetilde{\bf T}$, denoted by $\tau\sigma$, is obtained by identifying the roots of $\tau$ and $\sigma$ in the disjoint union $\tau\sqcup\sigma$, whose decorations are inherited from those of $\tau$ and $\sigma$ except that the $\frak{n}$-decoration at the root of $\tau\sigma$ is the sum of $\frak{n}$-decorations at the roots of $\tau$ and $\sigma$. We extend that tree product linearly to $\widetilde{T}$, which turns it into an algebra. 

For each $(\tau,{\frak{t}},{\frak{n}},{\frak{e}})\in\widetilde{\bf T}$ we define $\widetilde{\Delta}(\tau,{\frak{t}},{\frak{n}},{\frak{e}})$ as the infinite sum
\begin{equation} \label{*eq:graphicalcoprod}
\sum_\sigma\sum_{{\frak{n}}_\sigma,{\frak{e}}_{\partial\sigma}}\frac1{{\frak{e}}_{\partial\sigma}!}\binom{{\frak{n}}}{{\frak{n}}_\sigma} \big(\sigma,{\frak{t}}\vert_{E_\sigma},{\frak{n}}_\sigma+\pi{\frak{e}}_{\partial\sigma},{\frak{e}}\vert_{E_\sigma}\big) 
\otimes
\Big(\tau/\sigma,{\frak{t}}\vert_{E_\tau\setminus E_\sigma},[{\frak{n}}-{\frak{n}}_\sigma]_\sigma,{\frak{e}}\vert_{E_\tau\setminus E_\sigma}+{\frak{e}}_{\partial\sigma}\Big),
\end{equation}
where
\begin{itemize}
	\item[--] $\sigma$ runs over all subtrees of $\tau$ which contain the root of $\tau$. Then the quotient tree $\tau/\sigma$ is obtained by identifying all nodes of $\sigma$ in the tree $\tau$. The edge set of $\tau/\sigma$ is $E_\tau\setminus E_\sigma$.   \vspace{0.1cm}

	\item[--] ${\frak{n}}_\sigma$ runes over all maps $N_\sigma\to\bbN^{1+d}$ such that ${\frak{n}}_\sigma(v)\le{\frak{n}}(v)$ for any $v\in N_\sigma$. The map 
	$$
	[{\frak{n}}-{\frak{n}}_\sigma]_\sigma:N_{\tau/\sigma}\to\bbN^{1+d}
	$$ 
	is defined by $[{\frak{n}}-{\frak{n}}_\sigma]_\sigma(v)={\frak{n}}(v)$ for non-root $v\in N_{\tau/\sigma}$ and by $[{\frak{n}}-{\frak{n}}_\sigma]_\sigma(\varrho)=\sum_{v\in N_\sigma}({\frak{n}}(v)-{\frak{n}}_\sigma(v))$ for the root $\varrho$ of $\tau/\sigma$. Moreover one sets 
	$$
	\binom{{\frak{n}}}{{\frak{n}}_\sigma}\defeq\prod_{v\in N_\sigma}\binom{{\frak{n}}(v)}{{\frak{n}}_\sigma(v)}.
	$$ 

	\item[--] $\partial\sigma$ is the boundary of $\sigma$, that is, the set of all edges $(u,v)\in E_\tau$ such that $u\in N_\sigma$ and $v\in N_\tau\setminus N_\sigma$. ${\frak{e}}_{\partial\sigma}$ runs over all maps $\partial\sigma\to\bbN^{1+d}$.  For any $u\in N_\sigma$, $\pi{\frak{e}}_{\partial\sigma}(u)\defeq\sum_{(u,v)\in\partial\sigma}{\frak{e}}_{\partial\sigma}((u,v))$. Moreover one sets 
	$$
	{\frak{e}}_{\partial\sigma}!\defeq\prod_{e\in\partial\sigma}{\frak{e}}_{\partial\sigma}(e)!.
	$$
\end{itemize}
A proper definition of the infinite sum \eqref{*eq:graphicalcoprod} involves the notions of `bigraded space' and tensor product of bigraded spaces. See Section 2.3 of \cite{BHZ} for the details. Since we only consider below some truncations of $\widetilde{\Delta}$ into some finite sums, we do not touch the details here. It was shown in Proposition 3.23 of \cite{BHZ} that $\widetilde{T}$, equipped with the tree product and the coproduct $\widetilde\Delta$, is a Hopf algebra (as a bigraded space). 

The map $\widetilde\Delta$ admits a useful recursive characterization. A $4$-tuple $(\bullet,\varnothing,k,\varnothing)$ made up of a single node with decoration ${\frak{n}}(\bullet)=k\in \bbN^{1+d}$ is denoted by $X^k$. We use the notations ${\bf1}\defeq X^0$ and $X_j\defeq X^{e_j}$, where $e_j=(0,\dots,1,\dots,0)$ is the $j$-th vector of the canonical basis of $\bbN^{1+d}$. We write $\mcI_{k}^{\frak{l}}(\tau)$ for the planted tree obtained from $\tau\in\widetilde{\bf T}$ by grafting it to a new root with $\frak{n}$-decoration $0$, along an edge with label ${\frak{l}}\in{\frak{L}}$ and $\frak{e}$-decoration ${k}\in\bbN^{1+d}$; we consider $\mcI_k^{\frak{l}}:\widetilde{T}\to \widetilde{T}$ as a linear map. Then $\widetilde{\Delta}:\widetilde{T}\to \widetilde{T}\otimes \widetilde{T}$ is characterized by the recursive formulae
\begin{equation*}
\begin{aligned}
&\widetilde{\Delta}(\tau\sigma)=(\widetilde{\Delta}\tau)(\widetilde{\Delta}\sigma),  \qquad
\widetilde{\Delta}(X^k) = \sum_{l\le k} \binom{k}{l}X^{l}\otimes X^{k-l},   \\ 
&\widetilde{\Delta}(\mcI_k^{\frak{l}}(\tau)) = (\mcI_k^{\frak{l}}\otimes\textrm{Id})\widetilde{\Delta}\tau + \sum_{{l}\in\bbN^{1+d}}\frac{X^{l}}{{l}!}\otimes\mcI_{k+l}^{\frak{l}}(\tau)\qquad (\forall\,\tau,\sigma).
\end{aligned}
\end{equation*}
See Proposition 4.17 of \cite{BHZ} for a proof. Similarly to \cite{BH23}, we do not consider the trees that have some $\bbK$-type leaves. Because of this, we actually consider the quotient Hopf algebra
$$
T=\widetilde{T}/I,
$$
where $I$ is the Hopf ideal spanned by trees that have a $\bbK$-type leaf. We identify $T$ with the subspace spanned by the set $\bf T$ of all trees without $\bbK$-type leaves. Then the tree product $T\otimes T\to T$ and the coproduct 
$$
\Delta\defeq(\pi\otimes\pi)\circ\widetilde{\Delta}\circ\textit{inj}:T\to T\otimes T
$$
are well-defined, where $\textit{inj}:T\to\widetilde{T}$ is the natural injection map and $\pi:\widetilde{T}\to T$ is a canonical surjection. 

Fix $s_0\geq-\frac{2+d}2$, $\beta_0\in(0,2)$, and $r_0<-\frac{2+d}2-s_0$ such that 
$$
r_0\notin\bigg\{ 2\beta_0 + 2q_1 + \sum_{j=2}^{1+d} q_j \,;\, r,q_1,\dots,q_{1+d}\in\bbQ \bigg\}.
$$
(The lower bound of $s_0$ is not essential: it excludes the case $r_0>0$ where renormalization is no longer necessary. The parameter $\beta_0$ is chosen as an arbitrary number slightly smaller than $2$ because of the same technical reason as \cite{BH23}.) Recall the embeddings of function spaces over the $(1+d)$-dimensional space $\bbR^{1+d}$
$$
H^{-s_0,\mcQ}(w_b)\hookrightarrow B_{2,\infty}^{r_0+\frac{2+d}2,\mcQ}(w_a)\hookrightarrow B_{p,\infty}^{r_0+\frac{2+d}p,\mcQ}(w_a)\hookrightarrow C^{r_0,\mcQ}(w_a)=\Omega,
$$
for any $a\ge b\ge0$. For any $\varepsilon\ge0$ and $p\in[2,\infty]$, we define the degree map $r_{\varepsilon,p}:{\bf T}\to\bbR$ by

\begin{align*}
&r_{\varepsilon,p}(\bbOmega) = r_{\eps,p}(\bbH)\defeq r_0-\varepsilon,\qquad r_{\varepsilon,p}(\ubbH) \defeq r_0-\varepsilon+\frac{2+d}p, \qquad r_{\varepsilon,p}(\bbK) \defeq \beta_0,
\\
&r_{\varepsilon,p}(\tau_{\frak{e}}^{\frak{n}}) \defeq \sum_{v\in N_\tau}|{\frak{n}}(v)|_{\frak{s}} + \sum_{e\in E_\tau} \big(r_{\varepsilon,p}({\frak{t}}(e))-|{\frak{e}}(e)|_{\frak{s}}\big).
\end{align*}
The parameter $\epsilon$ is introduced for some technical purpose which will not be fundamental here. We refer to the work \cite{BH23} for more on this point. For each $\varepsilon\ge0$ and $p\in[2,\infty]$ we introduce the projection map $P^+_{\varepsilon,p}$ from $T$ to the subalgebra $T^+_{\varepsilon,p}$ spanned by the symbols $X^k\prod_{i=1}^n\mcI_{k_i}^{{\frak{l}}_i}(\tau_i)$ with $n\in\bbN$, $k,k_i\in\bbN^{1+d}$, ${\frak{l}}_i\in{\frak{l}}$ and $\tau_i\in{\bf T}$ such that $r_{\varepsilon,p}(\tau_i)+r_{\varepsilon,p}({\frak{l}}_i)>|k_i|_{\frak{s}}$ for each $i$. We define
\begin{equation*} \begin{split} 
&\Delta_{\varepsilon,p} \defeq(\textrm{Id}\otimes \mathop{P^+_{\varepsilon,p}})\Delta:T\to T\otimes T_{\varepsilon,p}^+,   \\
&\Delta^+_{\varepsilon,p} \defeq \big(P^+_{\varepsilon,p}\otimes P^+_{\varepsilon,p}\big)\Delta:T_{\varepsilon,p}^+\to T_{\varepsilon,p}^+\otimes T_{\varepsilon,p}^+.
\end{split} \end{equation*}
The ranges are algebraic tensor products since $P_{\varepsilon,p}^+$ restricts the choices of ${\frak{e}}_{\partial\sigma}$ in \eqref{*eq:graphicalcoprod}. By a similar argument to \cite{Hai14, BHZ} we can see that $T^+_{\varepsilon,p}$ is a Hopf algebra with coproduct $\Delta^+_{\varepsilon,p}$ and $T$ has a right comodule structure with coaction $\Delta_{\varepsilon,p}$. Denote by $S^+_{\varepsilon,p}$ the antipode of $(T^+_{\varepsilon,p},\Delta^+_{\varepsilon,p})$. 
The $p$-dependence of $\Delta_{\varepsilon,p}$ was described in an example (3.6) in \cite{BH23}.

\ssk

The linear space $T$ is too large for our purpose and we fix until the end of this section the subset ${\bf B}\subset\bfT$ of all trees $\tau\in\bfT$ with only $\bbK$ or $\bbOmega$ type edges that strongly conform to a given complete subcritical rule and such that $r_{0,\infty}(\tau) < C_0$ for a fixed number $C_0$. (Its choice is dictated by each equation when we apply the setting of regularity structures to the study of any given subcritical singular stochastic PDE.) Moreover, for each $(\tau,{\frak{t}},{\frak{n}},{\frak{e}})\in{\bf B}$, we assume that if $e=(u,v)\in {\frak{t}}^{-1}(\bbOmega)$ then $e$ is a leaf, ${\frak{n}}(v)={\frak{e}}(e)=0$, and there are no other edges $e'=(u,v')\in{\frak{t}}^{-1}(\bbOmega)$ with $v'\neq v$.

\medskip

\noindent \textbf{\S2. A regularity-integrability structure --} For each $(\tau,{\frak{t}},{\frak{n}},{\frak{e}})\in{\bf B}$ and disjoint subsets $V,\underline{V}$ of $\subset{\frak{t}}^{-1}(\bbOmega)$, we denote by $D_{V,\underline{V}}(\tau,{\frak{t}},{\frak{n}},{\frak{e}})$ the decorated tree $(\tau,{\frak{t}}',{\frak{n}},{\frak{e}})$ with the modified node decoration ${\frak{t}}'$ defined by
$$
{\frak{t}}'(e)=\begin{cases}
\bbH&(e\in V),\\
\ubbH&(e\in\underline{V}),\\
{\frak{t}}(e)&(\text{\rm otherwise}).
\end{cases}
$$ 
We set as well
\begin{align*}
\dot{\bf B} &\defeq \Big\{D_{V,\varnothing}\tau\ ;\, \tau\in{\bf B},\ \varnothing\neq V\subset{\frak{t}}^{-1}(\bbOmega) \Big\},   \\
\underline{\dot{\bf B}} &\defeq \Big\{ D_{V,\{v\}}\tau\ ;\, \tau\in{\bf B},\ V\subset{\frak{t}}^{-1}(\bbOmega),\ v\in {\frak{t}}^{-1}(\bbOmega)\setminus V \Big\},   \\
\widetilde{\bf B}&\defeq{\bf B}\cup\dot{\bf B}\cup\underline{\dot{\bf B}},
\end{align*}
and define the linear spaces
$$
V\defeq\spa({\bf B}),\qquad
U\defeq\spa({\bf B}\cup\dot{\bf B}),\qquad
W\defeq\spa(\widetilde{\bf B}).
$$
We set
\begin{align*}
    i_p(\tau) = 
    \begin{cases}
        \infty & \text{if\ } \tau \in \bfB \cup \dot{\bfB}   \\
        p & \text{if\ } \tau \in \underline{\dot{\bf B}}.
    \end{cases}
\end{align*}
For each $(\varepsilon,p)\in[0,\infty)\times[2,\infty]$ we define $V_{\varepsilon}^+$, $U_{\eps}^+$, and $W_{\eps,p}^+$ as the subalgebras of $T_{\varepsilon,p}^+$ generated by the families of symbols 
\begin{align*}
{\bf V}_{\varepsilon}^+&\defeq\{X_j\}_{j=1}^{1+d}\cup\{\mcI_k^{\bbK}(\tau)\}_{k\in\bbN^{1+d},\, \tau\in{\bf B}\setminus\{X^l\}_l,\, r_{\varepsilon,p}(\tau)+\beta_0>|k|_\mfs},\\
{\bf U}_{\varepsilon}^+&\defeq\{X_j\}_{j=1}^{1+d}\cup\{\mcI_k^{\bbK}(\tau)\}_{k\in\bbN^{1+d},\, \tau\in({\bf B}\cup\dot{\bf B})\setminus\{X^l\}_l,\, r_{\varepsilon,p}(\tau)+\beta_0>|k|_\mfs},\\
{\bf W}_{\varepsilon,p}^+&\defeq\{X_j\}_{j=1}^{1+d} \cup \big\{\mcI_k^{\ubbH}(\bullet)\big\}_{k\in\bbN^{1+d},\, r_{\varepsilon,p}(\ubbH)>|k|_\mfs} \cup \big\{\mcI_k^{\bbK}(\tau)\big\}_{k\in\bbN^{1+d},\, \tau\in\widetilde{\bf B}\setminus\{X^l\}_l,\, r_{\varepsilon,p}(\tau)+\beta_0>|k|_\mfs},
\end{align*}
respectively. We write ${\bf V}_\eps^+$ and ${\bf U}_\eps^+$ rather than ${\bf V}_{\eps,p}^+$ and ${\bf U}_{\eps,p}^+$ since $p$ has no influences on the trees in ${\bf B}\cup\dot{\bf B}$. (From its definition, the quantity $r_{\eps,p}(\tau)$ varies with $p$ only on trees with at least one $\ubbH$-type edge.)

\ssk

The proof of the following statement is omitted since it is almost identical to the proof of Lemma 6 in \cite{BH23}.

\ssk

\begin{lem} \label{LemActionCharacters}
The set $V_\eps^+$, resp. $U_\eps^+$ and $W_{\eps,p}^+$, is a sub-Hopf algebra of $T_{\eps,p}^+$, and $V$, resp. $U$ and $W$, is a right comodule over $V_\eps^+$, resp. $U_\eps^+$ and $W_{\eps,p}^+$, with coaction $\Delta_{\eps,p}$.
\end{lem}

\ssk

The group ${\sf G}_{\varepsilon,p}^+$ of characters on the Hopf algebra $(W_{\varepsilon,p}^+,\Delta_{\varepsilon,p}^+)$ has a representation in $GL(W)$ where ${\sf g}\in {\sf G}_{\varepsilon,p}^+$ is mapped to $(\id\otimes \mathop{\sf g})\Delta_{\varepsilon,p}$. Denote by ${\sf G}_{\varepsilon,p}$ the image group. Lemma \ref{LemActionCharacters} ensures that $V$ and $U$ are stable under the action of ${\sf G}_{\eps,p}$.

We finally define a grading on $Z\in\{V,U,W\}$ setting for any ${\bf a}\in\bbR\times[1,\infty]$ and ${\bf Z}\in\{{\bf B},\bfB\cup\dot{\bfB},\widetilde{\bfB}\}$
$$
Z_{\bf a}\defeq \spa\Big\{\tau \in {\bf Z}\, ;\, \big(r_{\varepsilon,p}(\tau),i_p(\tau)\big)={\bf a}\Big\}.
$$
With this grading, we define some regularity-integrability structures
\begin{align*}
\scrV_\eps &\defeq (A_{\eps},V,{\sf G}_{\eps,p}\vert_V),   \\
\scrU_\eps &\defeq (A_{\eps},U,{\sf G}_{\eps,p}\vert_U),   \\
\scrW_{\eps,p} &\defeq (A_{\eps,p},W,{\sf G}_{\eps,p}),
\end{align*}
where 
\[
A_\eps\defeq\big\{\big(r_{\varepsilon,\infty}(\tau),\infty\big)\, ;\, \tau\in{\bf B}\cup\dot{\bf B}\big\}\]
and 
\[
A_{\varepsilon,p} \defeq \big\{\big(r_{\varepsilon,p}(\tau),i_p(\tau)\big)\, ;\, \tau\in\widetilde{\bf B}\big\}.
\] 
These regularity-integrability structures are described below in their \textit{\textbf{concrete form}}
\begin{align*}
\scrV_\eps &= \big((V,\Delta_{\varepsilon,p}) , (V_\eps^+,\Delta_{\varepsilon,p}^+)\big),   \\
\scrU_\eps &= \big((U,\Delta_{\varepsilon,p}) , (U_\eps^+ , \Delta_{\varepsilon,p}^+)\big),   \\
\scrW_{\eps,p} &= \big((W,\Delta_{\varepsilon,p}) , (W_{\eps,p}^+,\Delta_{\varepsilon,p}^+)\big).
\end{align*}
Such a representation is useful when we want to make the generators of the structure group explicit.

\medskip

\noindent \textbf{\S3. Models on regularity-integrability structures --} We recall from \cite{BH23} some basic facts about models on regularity-integrability structures that are involved in the statement and the proof of Theorem \ref{ThmConstructionBPHZ} below. A reader familiar with regularity structures can skip this paragraph after looking at \eqref{EqDefnSizePiTau} and \eqref{EqDefnSizeGMu}. Denoting by $P_{\bf a} : T\to T_{\bf a}$ the canonical projection, we set with a slight abuse of notations
$$
\|\tau\|_{\bf a} \defeq \|P_{\bf a}\tau\|_{\bf a}
$$
for any $\tau\in T$ and ${\bf a}\in A$. 

\ssk

{\it Models --} Assume we are given a pair of maps ${\sf M}=({\sf \Pi}, {\sf g})$ such that
$$
{\sf\Pi} : W \to C^{r_0,\mcQ}(w_a),   \qquad   {\sf g} : \bbR^{1+d} \to {\sf G}_{\varepsilon,p}^+
$$
and $\sf\Pi$ is continuous and linear. The map $\sf \Pi$ is called an \textit{\textbf{interpretation map}}. For any $\tau\in W$ and $\mu\in W_{\varepsilon,p}^+$ set
\begin{equation*} \begin{split}
{\sf\Pi}^{\varepsilon,p}_x(\tau) &\defeq \big({\sf \Pi}\otimes {\sf g}_x^{-1}\big) \Delta_{\varepsilon,p}(\tau),   \\
{\sf g}_{yx}^{\varepsilon,p}(\mu) &\defeq \big({\sf g}_{y}\otimes {\sf g}_x^{-1}\big) \Delta^+_{\varepsilon,p}(\mu),
\end{split} \end{equation*}
where ${\sf g}_x^{-1}$ is an inverse of ${\sf g}_x$ in the group ${\sf G}_{\eps,p}^+$. A pair of maps ${\sf M}=({\sf \Pi}, {\sf g})$ as above is called a \textbf{\textit{model on $\scrW_{\varepsilon,p}$}} (with weight $w_a$) if
\begin{align} \label{EqDefnSizePiTau}
\|{\sf\Pi} \res \tau\|_{\epsilon,p;w_a} \defeq \sup_{0<t\le1} t^{-r_{\varepsilon,p}(\tau)/4} \big\| \mcQ_t(x, {\sf\Pi}^{\varepsilon,p}_x(\tau)) \big\|_{L_x^{i_p(\tau)}(w_a)}<\infty
\end{align}
for any $\tau\in\widetilde{\bf B}$ and
\begin{align} \label{EqDefnSizeGMu}
\|{\sf g} \res \mu\|_{\epsilon,p;w_a} \defeq \|{\sf g}_x(\mu)\|_{L_x^{i_p(\mu)}(w_a)} + \sup_{y\in\bbR^{1+d}\setminus\{0\}}\Bigg(w_a(y)\frac{\big\|{\sf g}^{\varepsilon,p}_{(x+y)x}(\mu)\big\|_{L_x^{i_p(\mu)}(w_a)}}{\|y\|_{\frak{s}}^{r_{\varepsilon,p}(\mu)}}\Bigg)<\infty
\end{align}
for any $\mu\in{\bf W}_{\varepsilon,p}^+$. 
For any models ${\sf M}=({\sf\Pi},{\sf g})$, we set
$$
\|{\sf M}\|_{\textbf{\textsf{M}}(\scrW_{\varepsilon,p})_{w_a}} \defeq \max_{\tau\in\widetilde{\bf B}} \|{\sf\Pi}\res \tau\|_{\epsilon,p;w_a}  +  \max_{\mu\in{\bf W}_{\varepsilon,p}^+} \|{\sf g} \res \mu\|_{\epsilon,p;w_a}.
$$
Moreover, we define a metric $\|{\sf M}_1 \res {\sf M}_2\|_{\textbf{\textsf{M}}(\scrW_{\varepsilon,p})_{w_a}}$ between two models ${\sf M}_1$ and ${\sf M}_2$ by
$$
\|{\sf M}_1 \res {\sf M}_2\|_{\textbf{\textsf{M}}(\scrW_{\varepsilon,p})_{w_a}} \defeq \max_{\tau\in\widetilde{\bf B}} \|{\sf\Pi}_1, {\sf\Pi}_2 \res \tau\|_{\epsilon,p;w_a}  +  \max_{\mu\in{\bf W}_{\varepsilon,p}^+} \|{\sf g}_1 , {\sf g}_2 \res \mu\|_{\epsilon,p;w_a}.
$$
The quantities in the right hand side are defined in the same way as \eqref{EqDefnSizePiTau} and \eqref {EqDefnSizeGMu} but with $({\sf\Pi}_1)_x^{\varepsilon,p}(\tau) - ({\sf\Pi}_2)_x^{\varepsilon,p}(\tau)$, $({\sf g}_1)_{x}(\mu) - ({\sf g}_2)_{x}(\mu)$ and $({\sf g}_1)_{yx}^{\varepsilon,p}(\mu) - ({\sf g}_2)_{yx}^{\varepsilon,p}(\mu)$ in place of ${\sf\Pi}^{\varepsilon,p}_x(\tau)$, ${\sf g}_x(\mu)$ and ${\sf g}^{\varepsilon,p}_{yx}(\mu)$, respectively. Finally, we define the space $\textbf{\textsf{M}}(\scrW_{\varepsilon,p})_{w_a}$ as the completion of collections of all models $\sf M=(\Pi,g)$ such that $\|{\sf M}\|_{\textbf{\textsf{M}}(\scrW_{\varepsilon,p})_{w_a}}$ is finite and ${\sf\Pi}\tau$ and ${\sf g}(\mu)$ are smooth for any $\tau$ and $\mu$, under the metric $\|\cdot \res \cdot\|_{\textbf{\textsf{M}}(\scrW_{\varepsilon,p})_{w_a}}$. This definition ensures that $\textbf{\textsf{M}}(\scrW_{\varepsilon,p})_{w_a}$ is a Polish space.

\ssk

$\bullet$ It will also be useful for later purposes to introduce a \textbf{\textit{homogeneous norm}} on the space of models on $\scrU_\eps$. Denote by $n(\tau)$ the number of edges of $\bbOmega$-type or $\bbH$-type in an arbitrary symbol $\tau$.
Set
\begin{align*}
\|{\sf\Pi} \res \tau\|_{\epsilon;w_a,\textrm{hom}} \defeq \left( \sup_{0<t\le1} t^{-r_{\varepsilon}(\tau)/4} \big\| \mcQ_t(x, {\sf\Pi}^{\varepsilon}_x(\tau)) \big\|_{L_x^{\infty}(w_a)} \right)^{1/n(\tau)},
\end{align*}
for all $\tau \in \bfB \cup\dot{\bfB}\setminus\{X^k\}_k$, and 
\begin{align*} 
\|{\sf g} \res \mu\|_{\epsilon;w_a,\textrm{hom}}
&\defeq 
\left ( \|{\sf g}_x(\mu)\|_{L_x^{\infty}(w_a)} + \sup_{y\in\bbR^{1+d}\setminus\{0\}}\Bigg(w_a(y) \, \frac{\big\|{\sf g}^{\varepsilon}_{(x+y)x}(\mu)\big\|_{L_x^{\infty}(w_a)}}{\|y\|_{\frak s}^{r_{\varepsilon}(\mu)}}\Bigg) \right)^{1/n(\mu)}
\end{align*}
for all $\mu\in{\bf U}_{\varepsilon}^+\setminus\{X^k\}_k$. Note that we can drop off the letter $p$ from these notations since $p$ has no influences for objects defined on $\scrU_\eps$.
We set
$$
\|{\sf M}\|_{\textbf{\textsf{M}}_{\textrm{hom}}(\scrU_{\varepsilon})_{w_a}} \defeq 
\max_{\tau\in(\bfB \cup \dot{\bfB})\setminus\{X^k\}_k} 
\|{\sf\Pi}\res \tau\|_{\epsilon;w_a,\textrm{hom}}  +  \max_{\mu\in{\bf U}_{\varepsilon}^+\setminus\{X^k\}_k} \|{\sf g} \res \mu\|_{\epsilon;w_a,\textrm{hom}}
$$
for any ${\sf M}\in\textbf{\textsf{M}}(\scrU_{\varepsilon})_{w_a}$. We also define similarly a metric $\|{\sf M}_1 \res {\sf M}_2\|_{\textbf{\textsf{M}}_{\textrm{hom}}(\scrU_{\varepsilon})_{w_a}}$ on $\textbf{\textsf{M}}(\scrU_{\varepsilon})_{w_a}$ and a metric $\|{\sf M}_1 \res {\sf M}_2\|_{\textbf{\textsf{M}}_{\textrm{hom}}(\scrV_{\varepsilon})_{w_a}}$ on $\textbf{\textsf{M}}(\scrV_{\varepsilon})_{w_a}$. The metrics $\|\cdot \res \cdot\|_{\textbf{\textsf{M}}(\scrV_{\varepsilon})_{w_a}}$ and $\|\cdot \res \cdot\|_{\textbf{\textsf{M}}_{\textrm{hom}}(\scrV_{\varepsilon})_{w_a}}$ generate the same topology on $\textbf{\textsf{M}}(\scrV_{\varepsilon})_{w_a}$. 

\ssk

$\bullet$ For the sub-concrete regularity-integrability structure  of the form
$$
\scrW_{\varepsilon,p}'=\big((Z,\Delta_{\varepsilon,p}),(Z^+,\Delta_{\varepsilon,p}^+)\big)
$$
where $Z=\spa({\bf C})$, for some subset ${\bf C}\subset\widetilde{\bf B}$, and $Z^+$ is a subalgebra of $W_{\varepsilon,p}^+$ generated by some subset ${\bf C}^+$ of ${\bf W}_{\varepsilon,p}^+$, it is useful to define the restricted quantity
\begin{equation*} 
\begin{aligned}
\|{\sf M}\|_{\textbf{\textsf{M}}(\scrW_{\varepsilon,p}')_{w_a}} \defeq 
\|{\sf\Pi}^{\varepsilon,p} \res {\bf C}\|_{\eps,p;w_a}  + \|{\sf g}^{\varepsilon,p} \res {\bf C}^+\|_{\eps,p;w_a} \defeq \max_{\tau\in\bf C} \|{\sf\Pi}^{\varepsilon,p} \res \tau\|_{\eps,p;w_a}
+\max_{\mu\in{\bf C}^+} \|{\sf g}^{\varepsilon,p} \res \mu\|_{\eps,p;w_a}.
\end{aligned}
\end{equation*}
In particular the restriction of any ${\sf M}\in\textbf{\textsf{M}}(\scrW_{\varepsilon,p})_{w_a}$ to $\scrV_\eps$ is a model in the usual sense of \cite{Hai14}. 

\ssk

{\it A family of renormalized models --} We choose $\bf B$ so that all the edges of an arbitrary $\tau\in\widetilde{\bf B}$ with labels other than $\bbK$ can be contracted into some labeled nodes. We write
$$
\ocircle=\mcI_0^{\bbOmega}(\bullet),\qquad
\odot=\mcI_0^{\bbH}(\bullet),\qquad
\underline{\odot}=\mcI_0^{\ubbH}(\bullet).
$$
Recall from \eqref{EqQSpacetimeOperator} the definition of the operator $\mcK$ acting on functions over $\bbR^{1+d}$. Denote by $C_\star^\infty$ be the set of all smooth functions over $\bbR^{1+d}$ whose derivatives of all orders are in the class $\bigcup_{a>0}L^\infty(w_a)$. An interpretation map $\sf \Pi$ is said to be \textit{\textbf{$\mcK$-admissible}} if it satisfies
$$
({\sf\Pi}X^k)(x)=x^k,\qquad{\sf\Pi}(\mcI_k^\bbK\tau) = \partial^k\mcK(\cdot,{\sf\Pi}\tau)
$$
for all $k\in\bbN^{1+d}$ and $\tau$. For any $\zeta,j\in C_\star^\infty$ there is a unique {\it multiplicative} $\mcK$-admissible map ${\sf\Pi}^{\zeta,j}:W\to C_\star^\infty$ such that
$$
{\sf\Pi}^{\zeta,j}(\ocircle) = \zeta,  \qquad  {\sf \Pi}^{\zeta,j}(\odot) = {\sf \Pi}^{\zeta,j}(\underline{\odot}) = j.
$$
We call a model ${\sf M}^{\zeta,j;\varepsilon,p} = ({\sf\Pi}^{\zeta,j}, {\sf g}^{\zeta,j;\varepsilon,p})$ on $\scrW_{\eps,p}$ (with weight $w_a$ for any large $a>0$) defined from the following recursive definition of $({\sf g}_x^{\zeta,j;\varepsilon,p})^{-1}$ the \textbf{\textit{naive model}} associated to $(\zeta,j)$

\begin{equation} \label{EqConstructionGfromPi}
\begin{aligned}
\big({\sf g}_x^{\zeta,j;\varepsilon,p}\big)^{-1}(X^k) &= (-x)^k,   \\
\big({\sf g}_x^{\zeta,j;\varepsilon,p}\big)^{-1}(\tau\sigma) &= ({\sf g}_x^{\zeta,j;\varepsilon,p})^{-1}(\tau)({\sf g}_x^{\zeta,j;\varepsilon,p})^{-1}(\sigma),   \\
\big({\sf g}_x^{\zeta,j;\varepsilon,p}\big)^{-1}(\mcI_k^{\frak{l}}(\sigma)) &= - \sum_{l\in\bbN^{1+d}} \, \frac{(-x)^{l}}{{l}!} {\sf f}_x^{\zeta,j;\varepsilon,p}(\mcI_{k+l}^{\frak{l}}(\sigma)),   \\
{\sf f}_x^{\zeta,j;\varepsilon,p}(\tau)
&={\bf1}_{r_{\varepsilon,p}(\tau)>0}\times
\left\{
\begin{aligned}
&\partial^k j(x)
&&(\tau=\mcI_k^{\ubbH}(\bullet)),   \\
&\partial^k\mcK(x, {\sf\Pi}_x^{\zeta,j;\varepsilon,p}\sigma)
&&(\tau=\mcI_k^\bbK(\sigma)).
\end{aligned}
\right.
\end{aligned}
\end{equation}
Denote by ${\bf B}_-$ the set of all $\tau\in{\bf B}$ such that $r_{0,\infty}(\tau)<0$ and $\tau\notin\{\ocircle\}\cup\{\mcI_k^{\bbK}(\sigma)\}_{k\in\bbN^{1+d},\,\sigma\in{\bf B}}$. We extend any map $\chi:{\bf B}_-\to\bbR$ into a linear map $\chi : W\to\bbR$ by setting $\chi(\tau)=0$ for $\tau\in\widetilde{\bf B}\setminus{\bf B}_-$, and define the linear map $R_\chi : W\to W$ from the formula
\begin{equation} \label{DefnRChi}
R_\chi(\tau) \defeq \tau + (\chi\otimes\id)\Delta\tau\qquad(\tau\in\widetilde{\bf B}).
\end{equation}
Given such a map $\chi$, we define $\widehat{M}^\chi:W\to W$ as the unique linear map fixing $X^k$, $\ocircle$, $\odot$, and $\underline{\odot}$, such that $\widehat{M}^\chi$ is {\it multiplicative} and satisfies
$$
\widehat{M}^\chi(\mcI_k^\bbK\tau) = \mcI_k^\bbK\big(\widehat{M}^\chi (R_\chi\tau)\big)
$$
for all $k$ and $\tau$. We also define the linear map
$$
M^\chi \defeq \widehat{M}^\chi R_\chi.
$$
Then for each naive model ${\sf M}^{\zeta,j;\varepsilon,p}$, we define the model
$$
{\sf M}^{\zeta,j,\chi;\varepsilon,p} \defeq \big( {\sf \Pi}^{\zeta,j,\chi}, {\sf g}^{\zeta,j,\chi;\varepsilon,p} \big)
$$
with 
\[
{\sf\Pi}^{\zeta,j,\chi} \defeq {\sf\Pi}^{\zeta,j}M^\chi
\]
and ${\sf g}^{\zeta,j,\chi;\varepsilon,p}$ defined by the version of the recursive rules \eqref{EqConstructionGfromPi} with $({\sf \Pi}^{\zeta,j,\chi}, {\sf g}^{\zeta,j,\chi;\varepsilon,p})$ in place of $({\sf \Pi}^{\zeta,j}, {\sf g}^{\zeta,j;\varepsilon,p})$. See Section 3 of Bruned's work \cite{Bru18} for the well-defined character of ${\sf M}^{\zeta,j,\chi;\varepsilon,p}$ as a model.

\ssk

\subsection{A BPHZ convergence result on the regularity-integrability structure $\scrW_{\eps,p}$. \hspace{0.1cm}}
\label{SubsectionBPHZConvergence}

We can now state and prove the BPHZ convergence theorem that constructs the BPHZ model which we will study in the next section. The statement and its proof are similar to the main result of \cite{BH23}; we will only describe below the structure of the argument, emphasizing the points where some modifications are needed.

\ssk

For a family of compactly supported smooth functions $\varrho_n\in C^\infty(\bbR^{1+d})$ that converge weakly to a Dirac mass at $0$ as $n\in\bbN$ goes to $\infty$, denote by $*$ the convolution operator and define a random variable $\xi_n$ on $\Omega$ setting for $\omega\in \Omega$ and $h\in H$
$$
\xi_n(\omega) \defeq \varrho_n*\omega \in\Omega, \qquad h_n  \defeq \varrho_n*h.
$$
We define an $\textbf{\textsf{M}}(\scrW_{\varepsilon,p})_{w_a}$-valued random variable setting for any deterministic map $\chi_n:\bfB_-\to\bbR$
$$
{\sf M}^{\xi_n,h_n,\chi_n;\varepsilon,p}(\omega)\defeq{\sf M}^{\xi_n(\omega),h_n,\chi_n;\varepsilon,p}.
$$ 
The main result of \cite{BHZ} yields that, for any $\bbP$ such that the law of $\xi_n$ is centered, transition invariant, has moments of all orders, we can find a unique $c_n:\bfB_-\to\bbR$ such that
\begin{equation*} 
\bbE\big[({\sf\Pi}^{\xi_n,h_n,\chi_n;0,\infty}\tau)(x)\big]=0
\end{equation*}
for all $\tau\in{\bf B}$ of non-positive $r_{0,\infty}$-degree. We call the renormalized model ${\sf M}^{\xi_n,h_n,\chi_n;\varepsilon,p}$ on the regularity-integrability structure $\scrW_{\eps,p}$ associated with this special choice of map $\chi_n$ the BPHZ model, and use for it the shorthand notation 
$$
{\sf M}^{n,h_n;\eps,p} = {\sf M}^{\xi_n,h_n,\chi_n;\varepsilon,p}.
$$
The following notation is involved in the next statement: ${\sf J}_p \defeq \big\{ \varepsilon\ge0\,;\,\mu\in{\bf W}_{0,2}^+,\ r_{\varepsilon,p}(\mu)=0 \big\}$, for an arbitrary $p\in[2,\infty]$. 

\ssk

\begin{thm} \label{ThmConstructionBPHZ}
Pick $s_0\geq-\frac{2+d}2$ and $r_0<-\frac{2+d}2-s_0$.
\begin{itemize}
	\item[--] Let $\bbP$ be a Borel probability measure on $\Omega=C^{r_0,\mcQ}(w_a)$, for some $a>0$, that satisfies the spectral gap inequality \eqref{EqSG} with $H = H^{-s_0,\mcQ}(w_b)$ for some $b\ge0$. 
	\item[--] Assume that 
	\begin{equation}\label{*BoundedVariance}
	r_{0,\infty}(\tau) > - \frac{2+d}{2}
	\end{equation}
	for all $\tau\in\bfB\setminus\{\ocircle\}$.
\end{itemize}	
Then there exist constants $\eps_0>0$ and $\overline{a}>0$ depending only on $b$, $d$, and ${\bf B}$ such that for any $a>\overline{a}$, $p\in[2,\infty], \varepsilon\in(0,\varepsilon_0)\setminus {\sf J}_p$ and $q\in[1,\infty)$, one has
\begin{equation*} 
\sup_{n\in\bbN}\bbE\Bigg[\sup_{\Vert h\Vert_H\leq 1} \Vert {\sf M}^{n,h_n;\varepsilon,p} \Vert_{\textbf{\textsf{M}}(\scrW_{\varepsilon,p})_{w_a}}^q\Bigg]<\infty.
\end{equation*}
and there exists a $\textbf{\textsf{M}}(\scrW_{\varepsilon,p})_{w_a}$-valued random variable ${\sf M}^{\infty,h;\eps,p}$ such that
$$
\lim_{n\to\infty}\bbE\bigg[\sup_{\Vert h\Vert_H\leq 1} \Vert{\sf M}^{n,h_n;\varepsilon,p} \res {\sf M}^{\infty,h;\varepsilon,p}\Vert_{\textbf{\textsf{M}}(\scrW_{\varepsilon,p})_{w_a}}^q\bigg]=0.
$$
The limit model ${\sf M}^{\infty,h;\eps,p}$ does not depend on the approximation of the identity $(\varrho_n)_{n\geq 1}$.
\end{thm}

\ssk

The use of the spectral gap assumption in the proof of convergence of the BPHZ random models was pioneered by Linares, Otto, Tempelmayr \& Tsatsoulis in \cite{LOTT}. We note that Theorem \ref{ThmConstructionBPHZ} is valid for spatially periodic noises $\xi$ on $\bbR\times\bbT^d$ satisfying Assumption A with $H=H^{-s_0}(\bbR\times\bbT^d)$, a usual flat $L^2$-Sobolev space with parabolic scaling. Indeed, since $H^{-s_0}(\bbR\times\bbT^d)\hookrightarrow H^{-s_0,\mcQ}(w_b)$ holds for any $b>0$, $\xi$ satisfies the spectral gap inequality \eqref{EqSG} with $H=H^{-s_0,\mcQ}(w_b)$. Although $H$ is of the form $H^{-s_0,\mcQ}(w_0)$ in \cite{BH23}, the same proof works for the weighted $H$ with only a slight modification of Lemma 15 of \cite{BH23}. When $\xi$ and $h$ are spatially periodic, then the model ${\sf M}^{n,h;\eps,p}$ is also spatially periodic in the sense that
$$
{\sf\Pi}_{x+e_i}^{n,h;\eps,p}(\tau)(\cdot+e_i)={\sf\Pi}_x^{n,h;\eps,p}(\tau)(\cdot),\qquad
{\sf g}_{(y+e_i)(x+e_i)}^{n,h;\eps,p}={\sf g}_{yx}^{n,h;\eps,p}
$$
for any $2\le i\le 1+d$. Recall that $e_i=(0,\dots,1,\dots,0)$ is the $i$-th vector of the canonical basis of $\bbN^{1+d}$.

\ssk

\begin{Dem}[of Theorem \ref{ThmConstructionBPHZ}]
The proof proceeds by induction, as illustrated in the picture below. We first decompose ${\bf B}=\{X^k\}_{|k|_\mfs<C_0}\cup{\bf B}_{\circ}$, where all $\tau\in\bfB_\circ$ has at least one $\bbOmega$-type edge. (The constant $C_0$ is involved in the description of the set $\bf B$.) Let $n_\circ(\tau)$ denote the number of noise symbols $\ocircle$ in an arbitrary symbol $\tau$. We align the elements of ${\bf B}_{\circ}$ as
$$
{\bf B}_{\circ} = \big\{\tau_1\preceq\tau_2\preceq\cdots\preceq\tau_N\big\},
$$
where $\sigma\preceq\tau$ means that
\begin{align*}
\big(n_\circ(\sigma), |E_\sigma|, r_{0,\infty}(\sigma)\big) \leq \big(n_\circ(\tau), |E_\tau|, r_{0,\infty}(\tau)\big)
\end{align*}
in the lexicographical order. We have in particular $\tau_1=\ocircle$. For $0\leq i\leq N$ set
$$
{\bf B}_i \defeq \big\{\tau_j\,;\,j\le i\big\}\cup\{X^k\}_{|k|_\mfs<C_0}.
$$
The following statement is proved as Lemma 14 in \cite{BH23}.

\ssk

\begin{lem} \label{leminduction}
For $1\le i\le N$, set
$$
V_i \defeq \spa({\bf B}_i),\qquad
U_i \defeq \spa\big({\bf B}_{i-1}\cup\dot{\bf B}_i\big),\qquad
W_i\defeq \spa\Big({\bf B}_{i-1}\cup\dot{\bf B}_{i-1}\cup\underline{\dot{\bf B}_i}\Big).
$$
For each $(\varepsilon,p)\in[0,\eps_0)\times[2,\infty]$ we define $V_{i,\varepsilon}^+$, $U_{i,\eps}^+$, and $W_{i,\eps,p}^+$ as the subalgebras of $T_{\varepsilon,p}^+$ generated by the families of symbols 
\begin{align*}
{\bf V}_{i,\varepsilon}^+&\defeq\{X_j\}_{j=1}^{1+d}\cup\{\mcI_k^{\bbK}(\tau)\}_{k\in\bbN^{1+d},\, \tau\in{\bf B}_i\setminus\{X^l\}_l,\, r_{\varepsilon,p}(\tau)+\beta_0>|k|_\mfs},\\
{\bf U}_{i,\varepsilon}^+&\defeq\{X_j\}_{j=1}^{1+d}\cup\{\mcI_k^{\bbK}(\tau)\}_{k\in\bbN^{1+d},\, \tau\in({\bf B}_i\cup\dot{\bf B}_i)\setminus\{X^l\}_l,\, r_{\varepsilon,p}(\tau)+\beta_0>|k|_\mfs},\\
{\bf W}_{i,\varepsilon,p}^+&\defeq\{X_j\}_{j=1}^{1+d} \cup \big\{\mcI_k^{\ubbH}(\bullet)\big\}_{k\in\bbN^{1+d},\, r_{\varepsilon,p}(\ubbH)>|k|_\mfs} \cup \big\{\mcI_k^{\bbK}(\tau)\big\}_{k\in\bbN^{1+d},\, \tau\in\widetilde{\bf B}_i\setminus\{X^l\}_l,\, r_{\varepsilon,p}(\tau)+\beta_0>|k|_\mfs},
\end{align*}
respectively. 
Then each of
\begin{align*}
\scrV_{i,\varepsilon} &\defeq \Big((V_i,\Delta_{\varepsilon,\infty}) , \big(V_{i-1}^+,\Delta_{\varepsilon,\infty}^+\big)\Big),  \\
\scrU_{i,\varepsilon} &\defeq \Big((U_i,\Delta_{\varepsilon,\infty}) , \big(U_{i-1}^+,\Delta_{\varepsilon,\infty}^+\big)\Big),  \\
\scrW_{i,\varepsilon,p} &\defeq \Big((W_i,\Delta_{\varepsilon,p}) , \big(W_{i-1,\varepsilon,p}^+,\Delta_{\varepsilon,p}^+\big)\Big)
\end{align*}
is a sub-regularity-integrability structure of $\scrW_{\varepsilon,p}$. 
\end{lem}

\ssk

For any fixed $p\in[2,\infty]$ we write below $\textsf{\textbf{bd}}(\scrW,i,p)$ to mean that
\begin{equation*} 
\sup_{n\in\bbN}\bbE\bigg[\sup_{\Vert h\Vert_H\leq 1} \Vert {\sf M}^{n,h_n;\varepsilon,p} \Vert_{\textbf{\textsf{M}}(\scrW_{i,\varepsilon,p})_{w_a}}^q\bigg]<\infty
\end{equation*}
for any small $\varepsilon>0$ and large $a>2+d$ and $q\in[1,\infty)$. Similarly, we write $\textsf{\textbf{cv}}(\scrW,i,p)$ to mean that 
$$
\lim_{n,m\to\infty}\bbE\bigg[\sup_{\Vert h\Vert_H\leq 1} \Vert{\sf M}^{n,h_n;\varepsilon,p} \res {\sf M}^{m,h_m;\varepsilon,p}\Vert_{\textbf{\textsf{M}}(\scrW_{\varepsilon,p})_{w_a}}^q\bigg]=0
$$
for any small $\varepsilon>0$ and large $a>2+d$ and all $q\in[1,\infty)$. We write 
$$
\textsf{\textbf{bd}}(\scrW,i), \quad  \textrm{respectively}  \quad \textsf{\textbf{cv}}(\scrW,i),
$$ 
to mean that $\textsf{\textbf{bd}}(\scrW,i,p)$, resp. $\textsf{\textbf{cv}}(\scrW,i,p)$ holds for any $p\in[2,\infty]$.
Moreover, we also write $\bound(\scrV,i)$, $\bound(\scrU,i)$, $\converge(\scrV,i)$ and $\converge(\scrU,i)$ to mean the analogue statements with $\scrV_{i,\varepsilon}$ or $\scrU_{i,\eps}$ in place of $\scrW_{i,\varepsilon,p}$. These statements do not depend on $p$ so we do not record this parameter in the notation. The following diagram summarizes the order of the proof. The dashed arrows mean that the corresponding steps are probabilistic. The solid arrows are independent of the probability measure.   
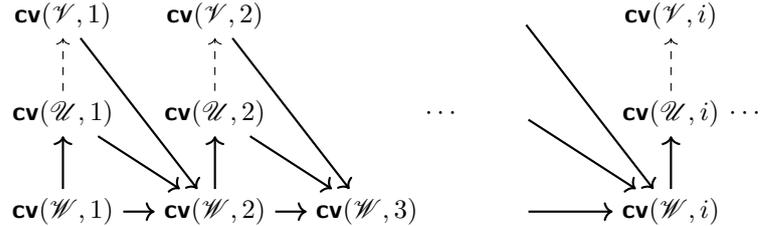
\vspace{0.15cm}
\begin{center}
\begin{tikzpicture}[auto]
\node (a1) at (0,0) {$\converge(\scrV,1)$}; 
\node (a2) at (2,0) {$\converge(\scrV,2)$}; 
\node (a3) at (4,0) {}; 
\node (a99) at (6,0) {}; 
\node (a100) at (8,0) {$\converge(\scrV,i)$}; 
\node (a101) at (9,0) {}; 
\node (b1) at (0,-1.3) {$\converge(\scrU,1)$}; 
\node (b2) at (2,-1.3) {$\converge(\scrU,2)$}; 
\node (b3) at (4,-1.3) {}; 
\node (b99) at (6,-1.3) {}; 
\node (b100) at (8,-1.3) {$\converge(\scrU,i)$}; 
\node (b101) at (9,-1.3) {}; 
\node (e1) at (0,-2.6) {$\converge(\scrW,1)$}; 
\node (e2) at (2,-2.6) {$\converge(\scrW,2)$}; 
\node (e3) at (4,-2.6) {$\converge(\scrW,3)$}; 
\node (e99) at (6,-2.6) {}; 
\node (e100) at (8,-2.6) {$\converge(\scrW,i)$}; 
\node (e101) at (9,-2.6) {}; 
\node (center) at ($0.5*(b3)+0.5*(b99)$) {$\cdots$};
\node (center2) at (b101) {$\cdots$};
\draw[->, thick] (e1) to (b1);
\draw[->, dashed] (b1) to (a1);
\draw[->, thick] (a1) to (e2);
\draw[->, thick] (e2) to (b2);
\draw[->, dashed] (b2) to (a2);
\draw[->, thick] (a2) to (e3);
\draw[->, thick] (a99) to (e100);
\draw[->, thick] (e100) to (b100);
\draw[->, dashed] (b100) to (a100);
\draw[->, thick] (e1) to (e2);
\draw[->, thick] (e2) to (e3);
\draw[->, thick] (e99) to (e100);
\draw[->, thick] (b1) to (e2);
\draw[->, thick] (b2) to (e3);
\draw[->, thick] (b99) to (e100);
\end{tikzpicture}
\end{center}

\ssk

The details of each step are as follows.

\ssk

\begin{itemize}
\setlength{\itemsep}{1mm}
\item[(0)] The initial case: $\converge(\scrW,1)$ is reduced to the fact that $h$ is an elements of $H = H^{-s_0,\mcQ}(w_b)$. This step is identical to the corresponding step in Section 4.2 of \cite{BH23}.
\item[(1)] $\converge(\scrW,i)\to\converge(\scrU,i)$: This step is trivial because $\converge(\scrU,i)$ can be seen as a particular case of $\converge(\scrW,i,\infty)$. Indeed, by the construction of models, we have
$$
{\sf \Pi}_x^{n,h_n;\varepsilon,\infty}(\tau) = {\sf\Pi}_x^{n,h_n;\varepsilon,\infty}(T(\tau))
$$
for any $\tau\in\underline{\dot{\bf B}}$, where $T:\underline{\dot{\bf B}}\to\dot{\bf B}$ is the linear map that replaces $\marked$ by $\odot$.

	\item[(2)] $\converge(\scrU,i)\to\converge(\scrV,i)$: This step is identical to the corresponding step in Section 4.3 of \cite{BH23}. 
\begin{itemize}
\item If $r_{0,\infty}(\tau_i)>0$, we obtain the result from $\converge(\scrV,i-1)$, which is included in $\converge(\scrU,i)$. This is a consequence of the reconstruction theorem. See Lemma 17 of \cite{BH23}.
\item Otherwise, we obtain a stochastic estimate from $\converge(\scrU,i)$ by using the spectral gap inequality. See Lemma 19 of \cite{BH23}.
\end{itemize}

	\item[(3)] $\big\{ \converge(\scrW,i),\converge(\scrU,i),\converge(\scrV,i) \big\} \to \converge(\scrW,i+1)$: This step is split into three parts.
\begin{itemize}
\setlength{\itemsep}{1mm}

	\item[(3-1)] First, we obtain some $n$-uniform bounds and the convergence of the $\sf g$-parts
$$
{\sf g}^{n,h_n;\eps,p}\big(\mcI_k^{\bbK}(\tau)\big)
$$
for any $\tau\in\widetilde{\bfB}_i$ from the assumptions $\big\{ \converge(\scrW,i),\converge(\scrU,i),\converge(\scrV,i) \big\}$; a consequence of the multilevel Schauder estimate in regularity-integrability structures. See Lemma 21 of \cite{BH23}.

	\item[(3-2)] Then one shows some $n$-uniform bounds and the convergence of the $\sf \Pi$-parts
$$
{\sf \Pi}^{n,h_n;\eps,p}(\tau)
$$
for any $\tau\in\underline{\dot{\bfB}_{i+1}}$. It is important for that purpose that the second assumption of Theorem \ref{ThmConstructionBPHZ} implies that $r_{0,2}(\tau)>0$ for any $\tau\in\widetilde{\bfB}\setminus\{\ocircle,\odot,\marked\}$. Hence, if $p=2$, the result is obtained as an application of the Reconstruction Theorem in regularity-integrability structures. See Lemma 22 of \cite{BH23}.

	\item[(3-3)] To extend the result $\converge(\scrW,i+1,2)$ into $\converge(\scrW,i+1,p)$ for all $p\in[2,+\infty]$, we need a comparison formula between ${\sf\Pi}^{n,h_n;\eps,2}$ and ${\sf\Pi}^{n,h_n;\eps,p}$ -- see Lemma 8 of \cite{BH23}. The proof of this lemma is exactly the same as in \cite{BH23}, except that the role of $\bbH$ is replaced here by $\ubbH$.
Once that formula is established, the proof of $\converge(\scrW,i+1,p)$ follows the same steps as the proof of Lemma 24 in \cite{BH23}.
\end{itemize}
\end{itemize}
This concludes the description of the  structure of the proof of Theorem \ref{ThmConstructionBPHZ}.
\end{Dem}

\ssk

Despite some similar notations, the spaces $\scrW$ and $\scrV$ above are not the same as the corresponding spaces of \cite{BH23}. The above statement is strictly stronger than the corresponding statement in \cite{BH23}, and Theorem \ref{ThmConstructionBPHZ} cannot be obtained as a direct corollary of Theorem $10$ in \cite{BH23}.

\medskip

\section{Transportation cost inequalities for the BPHZ models}
\label{SectionTCIModels}

We discuss in this section the limit BPHZ models on $\scrV_\eps$ and $\scrU_\eps$, so we can drop off from the notations the letter $p$ in the sequel. When discussing the restriction of BPHZ models to $\scrV_\eps$, we also write ${\sf M}^{n;\eps}$ instead of ${\sf M}^{n,h_n;\eps}$ since this restriction is independent of $h$. By Theorem \ref{ThmConstructionBPHZ}, there exists a subsequence $\big({\sf M}^{n_k,h_{n_k};\eps}(\omega)\big)_{k\geq 1}$ such that
$$
\lim_{k\to\infty}\sup_{\Vert h\Vert_H\leq 1} \big\Vert {\sf M}^{n_k,h_{n_k};\varepsilon}(\omega) \res {\sf M}^{\infty,h;\varepsilon}(\omega) \big\Vert_{\textbf{\textsf{M}}(\scrU_\varepsilon)_{w_a}}=0
$$
almost surely. We denote by $\Omega_0$ the set of all $\omega\in\Omega$ for which this convergence holds. By the multilinearity for $h$, if $\omega\in\Omega_0$, then $\big({\sf M}^{n_k,h_{n_k};\eps}(\omega)\big)_{k\geq 1}$ converges for any $h\in H$.
We define
$$
\overline{\sf M}^{\infty,h;\eps}(\omega) =
\begin{cases}
\lim_{k\to\infty}{\sf M}^{n_k,h_{n_k};\eps}(\omega)&(\omega\in\Omega_0),\\
{\bf 0}&(\omega\notin\Omega_0).
\end{cases}
$$
Recall that $\bf 0$ is the zero model over $\scrU_\eps$. On the other hand, we denote by $\Omega_1$ the set of all $\omega\in\Omega$ for which the restrictions ${\sf M}^{n_k;\eps}(\omega)$ of ${\sf M}^{n_k,h_{n_k};\eps}(\omega)$ into a model on $\scrV_\varepsilon$ converge as $k$ goes to $\infty$.
We define
$$
\widetilde{\sf M}^{\infty;\eps}(\omega) =
\begin{cases}
\lim_{k\to\infty}{\sf M}^{n_k;\eps}(\omega)&(\omega\in\Omega_1),\\
{\bf 0}&(\omega\notin\Omega_1).
\end{cases}
$$
Note that $\Omega_0\subset\Omega_1$ and the restriction of $\overline{\sf M}^{\infty,h;\eps}(\omega)$ onto $\scrV_\eps$ coincides with $\widetilde{\sf M}^{\infty;\eps}(\omega)$ for any $\omega\in\Omega_0$. Set 
\[
m({\bf B}) \defeq \max \big\{n_\circ(\tau) \,;\, \tau \in \bfB\big\}.
\] 

\ssk

\begin{thm} \label{ThmHolderEstimateBPHZ}
One has $\Omega_0+h\subset\Omega_1$ for all $h\in H$, and
 \begin{align}\label{mainbound}
        \big\| {\widetilde{\sfM}}^{\infty;\eps}(\omega+h) : {\widetilde{\sfM}}^{\infty;\eps}(\omega) \big\|_{\textbf{\textsf{M}}_{\textrm{\emph{hom}}}(\scrV_{\varepsilon})_{w_a}} 
        &\leq \mathbb{L}(\omega) \max\big(\|h\|_H,\|h\|^{1/m({\bf B})}_H\big)
\end{align}
for all $\omega\in\Omega_0$ and $h\in H$, for a random variable $\mathbb{L}(\omega)$ proportional to 
\[
\sup_{\|h\|_H\le1} \big\| {\overline{\sfM}}^{\infty,h;\eps}(\omega) \big\|_{\textbf{\textsf{M}}_{\textrm{\emph{hom}}}(\scrU_{\varepsilon})_{w_a}}
\] 
on $\Omega_0$ and infinite elsewhere.
\end{thm}

\ssk

We note that the random variable $\mathbb{L}$ is in all the $L^p(\bbP)$ spaces for $1\leq p<\infty$, from Theorem \ref{ThmConstructionBPHZ}. Hence the random variable ${\widetilde{\sfM}}^{\infty;\eps}:\Omega\to\textbf{\textsf{M}}(\scrV_{\varepsilon})_{w_a}$ satisfies the inequality \eqref{EqConditionExtendedContractionPrinciple} for all $\omega_1\in\Omega$ and $\omega_2\in\Omega_0$. It follows from Proposition \ref{PropGeneralisedTCI} that the law of $\widetilde{\sfM}^{\infty;\varepsilon}$ on the space $\textbf{\textsf{M}}(\scrV_{\varepsilon})_{w_a}$ satisfies for any $1\leq \alpha<2$ an $(\ell_\alpha , c^\alpha_{\textbf{\textsf{M}}})$-transportation cost inequality \eqref{EqTCI} with $\ell_\alpha$ the convex and continuous inverse function of the non-negative increasing continuous concave function $(t\in [0,\infty))\mapsto a_\alpha^{-1} (t^{\alpha/(2r)} + t^{\alpha/2})$. 

\ssk

The inequality \eqref{mainbound} is obtained from the following algebraic relations. For the sake of generality, we consider a deterministic model ${\sf M}^{\zeta,j,\chi}$ defined by some inputs $\zeta,j\in C_\star^\infty$ and some renormalization map $\chi:\bfB_-\to\bbR$. We suppress from the notations the exponents $\varepsilon,p$ to simplify the notations here and in the remainder of this section; these exponents have some fixed values.

\ssk

\begin{prop} \label{prop:modelbinomial}
For all $x,y$, for any $\tau\in\bfB$ one has
\begin{align} \label{eq:modelbinomialexpansion}
{\sf\Pi}_x^{\zeta+j,0,\chi}(\tau) = \sum_{V\subset\ocircle_\tau}{\sf\Pi}_x^{\zeta,j,\chi}(D_V\tau),
\end{align}
and for any $\sigma\in {\bf V}^+$ one has
\begin{align} \label{prop:modelbinomial2}
{\sf g}_{yx}^{\zeta+j,0,\chi}(\sigma)=\sum_{V\subset\ocircle_\sigma}{\sf g}_{yx}^{\zeta,j,\chi}(D_V\sigma),
\end{align}
where $\ocircle_\tau\defeq\frak{t}^{-1}(\bbOmega)$ and $D_V\tau\defeq D_{V,\varnothing}\tau$.
\end{prop}

\ssk

\begin{Dem}[of Theorem \ref{ThmHolderEstimateBPHZ} from Proposition \ref{prop:modelbinomial}]
First note that, if $\omega\in\Omega_0$, the subsequence $\{{\sf M}^{\xi_{n_k}+h_{n_k},0,\chi_{n_k}}\}$ converges via the decompositions \eqref{eq:modelbinomialexpansion} and \eqref{prop:modelbinomial2}, so $\omega+h\in\Omega_1$.

As an application of \eqref{eq:modelbinomialexpansion} we have
\begin{align*}
\big\| {\sf\Pi}^{n;\eps}(\omega+h) &, {\sf\Pi}^{n;\eps}(\omega)\res\tau \big\|_{\epsilon,p;w_a,\text{hom}} \notag   \\
&=\left (\sup_{0<t\le1}t^{-r_{\eps}(\tau)/4} \big\| \mcQ_t\big(x,{\sf\Pi}_x^{\xi_n+h_n,0,\chi_n}(\tau) - {\sf\Pi}_x^{\xi_n,0,\chi_n}(\tau)\big) \big\|_{L_x^\infty(w_a)} \right)^{\frac1{n(\tau)}}   \notag   \\
& \le \sum_{\varnothing\neq V\subset\ocircle_\tau} \left(\sup_{0<t\le1}t^{-r_{\eps}(\tau)/4} \big\| \mcQ_t\big(x,{\sf\Pi}_x^{\xi_n,h_n,\chi_n}(D_V\tau) \big) \big\|_{L_x^\infty(w_a)}\right)^{\frac1{n(\tau)}}   \notag   \\
&\le\sum_{\varnothing\neq V\subset\ocircle_\tau}\bigg(\sup_{\|h\|_H\le1}\|{\sf\Pi}^{n,h_n;\eps}\res D_V\tau\|_{\epsilon,p;w_a,\text{hom}}\bigg) \|h\|_H^{\frac{|V|}{n(\tau)}}.
\end{align*}
In the last inequality we use the fact that ${\sf\Pi}^{n,h;\eps}(D_V\tau)$ is homogeneous with respect to $h$. 
Since $\frac1{m(\bfB)}\le\frac{|V|}{n_\circ(\tau)}\le1$, we have
\begin{align*}
\big\| {\sf\Pi}^{n;\eps}&(\omega+h) , {\sf\Pi}^{n;\eps}(\omega)\res\tau \big\|_{\epsilon,p;w_a,\text{hom}}   \\
& \lesssim \left(\sup_{\|h\|_H\leq 1}\|{\sf\Pi}^{n,h_n;\varepsilon}(\omega)\|_{\textbf{\textsf{M}}_{\text{hom}}(\scrU_{\varepsilon})_{w_a}} \right) \max\big( \|h\|_H,\|h\|_H^{1/m({\bf B})}\big).
\end{align*}
We obtain a similar estimate for $\big\| {\sf g}^{n;\eps}(\omega+h) , {\sf g}^{n;\eps}(\omega)\res\mu \big\|_{\epsilon,p;w_a,\text{hom}}$ starting from \eqref{prop:modelbinomial2}. Taking the (subsequential) limit, we have the inequality \eqref{mainbound}. 
\end{Dem}

\ssk

We need the following two lemmas for the proof of Proposition \ref{prop:modelbinomial}. Recall from \eqref{DefnRChi} the definition of the linear map $R_\chi$.

\ssk

\begin{lem} \label{lem:preparation}
One has $R_\chi(D_V\tau) = D_V(R_\chi(\tau))$ for all $\tau\in\bfB$ and $V\subset\ocircle_\tau$. In the right hand side, $D_V$ acts only on $\sigma\in\bfB$ such that $V\subset\ocircle_\sigma$. Otherwise $D_V\sigma=0$.
\end{lem}

\ssk

\begin{Dem}
Let us use the notation $\Delta\tau=\sum_{\tau_1,\tau_2} \delta(\tau_1,\tau_2) \tau_1\otimes\tau_2$ with a sum over decorated trees $\tau_1,\tau_2$ and with some coefficients $\delta(\tau_1,\tau_2)$. Since $\Delta$ acts on $\bbOmega$-type edges and $\bbH$-type edges in the same way we can write
$$
\Delta(D_V\tau) = \sum_{\tau_1,\tau_2} \delta(\tau_1,\tau_2) \big(D_{V\cap \ocircle_{\tau_1}}(\tau_1)\big)\otimes \big(D_{V\cap \ocircle_{\tau_2}}(\tau_2)\big).
$$
Since $\chi$ vanishes on all the trees with at least one $\bbH$-type edge we have the consequence by applying $\chi\otimes\textrm{Id}$ to the above formula. Indeed one has
\begin{align*}
(\chi\otimes\textrm{Id})\Delta(D_V\tau)
=\sum_{\tau_1,\tau_2;V\subset\ocircle_{\tau_2}} \delta(\tau_1,\tau_2) \chi(\tau_1)\otimes \big(D_V(\tau_2)\big)
=D_V(R_\chi(\tau)).
\end{align*}  
\end{Dem}

\ssk

The next formula appears in Section 5.5 of \cite{BH23}; we state it here with a proof for completeness. Recall from \eqref{EqConstructionGfromPi} the definition of the character ${\sf f}_x^{\zeta,j,\chi}$. We extend ${\sf g}_x^{\zeta,j,\chi}$ and ${\sf f}_x^{\zeta,j,\chi}$ linearly and multiplicatively by setting ${\sf g}_x^{\zeta,j,\chi}(\eta) = {\sf f}_x^{\zeta,j,\chi}(\eta) = 0$ for any planted trees $\eta$ such that $r_{\eps,p}(\eta)\le0$.

\ssk

\begin{lem} \label{LemTwo}
For any $\sigma\in\bfB$ and $k\in\bbN^{1+d}$ we have 
\begin{align}\label{eq:express_g_by_f}
{\sf g}_{yx}^{\zeta,j,\chi}(\mcI_k^\bbK\sigma) = ({\sf f}_y^{\zeta,j,\chi}\otimes{\sf g}_{yx}^{\zeta,j,\chi})(\mcI_k^\bbK\otimes\textrm{\emph{Id}})(\Delta\sigma) - \sum_{l\in\bbN^{1+d}}\frac{(y-x)^l}{l!} \, {\sf f}_x^{\zeta,j,\chi}(\mcI_{k+l}^\bbK\sigma).
\end{align}
\end{lem}

\ssk

\begin{Dem}
Let us write in this proof ${\sf g}_x={\sf g}_x^{\zeta,j,\chi}$ and ${\sf f}_x = {\sf f}_x^{\zeta,j,\chi}$. By applying the operator ${\sf g}_x\otimes{\sf g}_x^{-1}\otimes{\sf g}_x$ to the identity
\begin{align*}
\sum_{m\in\bbN^{1+d}}\frac{X^m}{m!}\otimes\Delta^+(\mcI_{k+m}^{\bbK}\sigma)
=\sum_{m\in\bbN^{1+d}}\frac{X^m}{m!}\otimes(\mcI_{k+m}^{\bbK}\otimes\textrm{Id})\Delta\sigma
+\sum_{m,l\in\bbN^{1+d}}\frac{X^m}{m!}\otimes\frac{X^l}{l!}\otimes\mcI_{k+l+m}^{\bbK}\sigma,
\end{align*}
and using the fact $({\sf g}_x^{-1}\otimes{\sf g}_x)\Delta^+(\mcI_{k+m}^{\bbK}\sigma)=0$ and the binomial theorem, we have the explicit representation of ${\sf g}_x$
\begin{align*}
{\sf g}_x(\mcI_k^{\bbK}\sigma) = ({\sf f}_x\otimes{\sf g}_x)(\mcI_k^{\bbK}\otimes\textrm{Id})\Delta\sigma.
\end{align*}
Then we have as a consequence
\begin{align*}
({\sf f}_y\otimes{\sf g}_{yx})(\mcI_k^\bbK\otimes\id)\Delta\sigma
&=({\sf f}_y\otimes{\sf g}_y\otimes{\sf g}_x^{-1})(\mcI_k^{\bbK}\otimes\textrm{Id}\otimes\textrm{Id})(\id\otimes\Delta^+)\Delta\sigma   \\
&=({\sf f}_y\otimes{\sf g}_y\otimes{\sf g}_x^{-1})(\mcI_k^{\bbK}\otimes\textrm{Id}\otimes\textrm{Id})(\Delta\otimes\id)\Delta\sigma   \\
&=({\sf g}_y\otimes{\sf g}_x^{-1})(\mcI_k^{\bbK}\otimes\textrm{Id})\Delta\sigma   \\
&=({\sf g}_y\otimes{\sf g}_x^{-1})\bigg(\Delta^+\mcI_k^{\bbK}\sigma - \sum_{l\in\bbN^{1+d}}\frac{X^l}{l!}\otimes \mcI_{k+l}^{\bbK}\sigma\bigg)   \\
&={\sf g}_{yx}(\mcI_k^{\bbK}\sigma)+\sum_{l,m\in\bbN^{1+d}}\frac{y^l}{l!}\frac{(-x)^m}{m!} \, {\sf f}_x(\mcI_{k+l+m}^{\bbK}\sigma)   \\
&={\sf g}_{yx}(\mcI_k^{\bbK}\sigma)+\sum_{b\in\bbN^{1+d}}\frac{(y-x)^b}{b!} \, {\sf f}_x(\mcI_{k+b}^{\bbK}\sigma),
\end{align*}
which concludes the proof of the lemma.
\end{Dem}

\ssk

\begin{Dem}[of Proposition \ref{prop:modelbinomial}]
Recall from Proposition 3.15 in Bruned's work \cite{Bru18} that we can factorize the renormalized interpretation map as
\begin{equation*}
{\sf\Pi}_x^{\zeta,j,\chi} = \widehat{\sf\Pi}_x^{\zeta,j,\chi}R_\chi,
\end{equation*}
where $\widehat{\sf\Pi}_x^{\zeta,j,\chi}$ is the linear and multiplicative map defined by $\widehat{\sf\Pi}_x^{\zeta,j,\chi}X^k = (\cdot-x)^k$ and
$$
\widehat{\sf\Pi}_x^{\zeta,j,\chi}(\mcI_k^\bbK\tau) = \partial^k\mcK(\cdot,{\sf\Pi}_x^{\zeta,j,\chi}\tau)-\sum_{l\in\bbN^{1+d},\,|l|_\mfs<r(\mcI_k^\bbK\tau)}\frac{(\cdot-x)^l}{l!}\partial^{k+l}\mcK\big(x,{\sf\Pi}_x^{\zeta,j,\chi}\tau\big).
$$
First we prove \eqref{eq:modelbinomialexpansion} and the auxiliary equality
\begin{equation}\label{eq:hatmodelbinomialexpansion}
\widehat{\sf\Pi}^{\zeta+j,0,\chi}(\tau) = \sum_{V\subset\ocircle_\tau}\widehat{\sf\Pi}^{\zeta,j,\chi}(D_V\tau)
\end{equation}
simultaneously. The proof is an induction on the size of $\tau$.
Both \eqref{eq:modelbinomialexpansion} and \eqref{eq:hatmodelbinomialexpansion} are obvious for the initial cases $\tau\in\{\ocircle,X^k\}$. 
Next let $\tau$ be a planted tree of the form $\mcI_k^\bbK(\sigma)$ with $\sigma\in{\bf B}$. If $\sigma$ satisfies \eqref{eq:modelbinomialexpansion} then, by the definition of the operator $\widehat{\sf\Pi}_x^{\zeta,j,\chi}$, we have
\begin{align*}
\widehat{\sf\Pi}_x^{\zeta+j,0,\chi}(\tau)
&=\partial^k\mcK\big(\cdot,{\sf\Pi}_x^{\zeta+j,0,\chi}(\sigma)\big) - \sum_{|l|_\mfs<r(\tau)}\frac{(\cdot-x)^l}{l!}\partial^{k+l}\mcK\big(x,{\sf\Pi}_x^{\zeta+j,0,\chi}(\sigma)\big)   \\
&=\sum_{V\subset\ocircle_\sigma}\bigg(\partial^k\mcK\big(\cdot,{\sf\Pi}_x^{\zeta,j,\chi}(D_V\sigma)\big) - \sum_{|l|_\mfs<r(\tau)}\frac{(\cdot-x)^l}{l!}\partial^{k+l}\mcK\big(x,{\sf\Pi}_x^{\zeta,j,\chi}(D_V\sigma)\big)\bigg)   \\
&=\sum_{V\subset\ocircle_\sigma}\widehat{\sf\Pi}_x^{\zeta,j,\chi}(\mcI_k^\bbK(D_V\sigma)). 
\end{align*}
In the last equality we use the fact that $D_V\sigma$ is a linear combination of decorated trees with degree $r(\tau)$. Since $\ocircle_\sigma=\ocircle_\tau$ and $\mcI_k^\bbK D_V\sigma=D_V\mcI_k^\bbK\sigma$, we obtain that $\tau$ satisfies \eqref{eq:hatmodelbinomialexpansion}. For any non-planted trees $\tau=\prod_{i=0}^N\eta_i$, we also obtain \eqref{eq:hatmodelbinomialexpansion} by the multiplicativity of $\widehat{\sf\Pi}_x^{\zeta,j,\chi}$. Finally, by using the result of Lemma \ref{lem:preparation}, we see that
\begin{align*}
{\sf\Pi}^{\zeta+j,0,\chi}(\tau)
&=
\widehat{\sf\Pi}^{\zeta+j,0,\chi}(R_\chi\tau)
=
\sum_{V\subset\ocircle_\tau}\widehat{\sf\Pi}^{\zeta,j,\chi}(D_V(R_\chi\tau))   \\
&=
\sum_{V\subset\ocircle_\tau}\widehat{\sf\Pi}^{\zeta,j,\chi}(R_\chi(D_V\tau))
=
\sum_{V\subset\ocircle_\tau}{\sf\Pi}^{\zeta,j,\chi}(D_V\tau).
\end{align*}
Finally, we prove \eqref{prop:modelbinomial2} by an induction on the size of trees. By multicativity, it is sufficient to consider the planted tree $\sigma=\mcI_k^\bbK(\eta)$. Note that ${\sf f}_x^{\zeta,j,\chi}$ satisfies from \eqref{eq:modelbinomialexpansion} the formula
$$
{\sf f}_x^{\zeta+j,0,\chi}(\tau)=\sum_{V\subset\ocircle_\tau}{\sf f}_x^{\zeta,j,\chi}(D_V\tau).
$$
By applying \eqref{eq:express_g_by_f} of Lemma \ref{LemTwo} to $\Delta(\eta) \eqdef \sum_{\eta_1,\eta_2} \delta(\eta_1,\eta_2) \eta_1\otimes\eta_2$, for some decorated trees $\eta_1,\eta_2$ and coefficients $\delta(\eta_1,\eta_2)$, we have

\begin{align*}
{\sf g}_{yx}^{\zeta+j,0,\chi}(\mcI_k^\bbK\eta) &= \sum_{\eta_1,\eta_2} \delta(\eta_1,\eta_2) {\sf f}_y^{\zeta+j,0,\chi}(\mcI_k^\bbK\eta_1) \, {\sf g}_{yx}^{\zeta+j,0,\chi}(\eta_2) - \sum_{l\in\bbN^{1+d}}\frac{(y-x)^l}{l!} \, {\sf f}_x^{\zeta+j,0,\chi}(\mcI_{k+l}^\bbK\eta)   \\
&= \sum_{\eta_1,\eta_2} \delta(\eta_1,\eta_2) \sum_{V_1\subset\ocircle_{\eta_1},V_2\subset\ocircle_{\eta_2}} {\sf f}_y^{\zeta,j,\chi}(\mcI_k^\bbK D_{V_1}\eta_1) \, {\sf g}_{yx}^{\zeta,j,\chi}(D_{V_2}\eta_2)   \\
&\quad - \sum_{l\in\bbN^{1+d}}\frac{(y-x)^l}{l!}\sum_{V\subset\ocircle_\eta} {\sf f}_x^{\zeta,j,\chi}\big(\mcI_{k+l}^\bbK(D_V\eta)\big)   \\
&= \sum_{V\subset\ocircle_\eta} \bigg( \sum_{\eta_1,\eta_2} \delta(\eta_1,\eta_2) \, {\sf f}_y^{\zeta,j,\chi}(\mcI_k^\bbK D_{V\cap\ocircle_{\eta_1}}\eta_1) \, {\sf g}_{yx}^{\zeta,j,\chi} (D_{V\cap\ocircle_{\eta_2}}\eta_2)   \\
&\quad - \sum_{l\in\bbN^{1+d}}\frac{(y-x)^l}{l!} \, {\sf f}_x^{\zeta,j,\chi}\big(\mcI_{k+l}^\bbK(D_V\eta)\big) \bigg)   \\
&= \sum_{V\subset\ocircle_\eta}\bigg( ({\sf f}_y^{\zeta,j,\chi}\otimes{\sf g}_{yx}^{\zeta,j,\chi})(\mcI_k^\bbK\otimes\id)\Delta(D_V\eta) - \sum_{l\in\bbN^{1+d}}\frac{(y-x)^l}{l!}{\sf f}_x^{\zeta,j,\chi}\big(\mcI_{k+l}^\bbK(D_V\eta)\big) \bigg)   \\
&= \sum_{V\subset\ocircle_\eta}{\sf g}_{yx}^{\zeta,j,\chi}(\mcI_k^\bbK\eta).
\end{align*}
This concludes the proof of Proposition \ref{prop:modelbinomial}.
\end{Dem}

\ssk

We note that Sch\"onbauder \cite{Sch23} also studied the extended BPHZ model including Malliavin derivatives of all orders, to show Malliavin differentiability of solutions to singular stochastic PDEs. In \cite{Sch23}, building up the convergence of the BPHZ model over $\scrV_\eps$ from \cite{ChandraHairer}, a lift map from $h\in H$ to the extended BPHZ model $\widetilde{\sfM}^{\infty,h;\eps}$ over $\scrU_\eps$ is constructed. In contrast, in our proof, BPHZ models over $\scrV_\eps$ and $\scrU_\eps$ are constructed simultaneously in a self-contained proof. Moreover, while only Gaussian noises are considered in \cite{Sch23}, our proof is valid for all noises satisfying the spectral gap inequality. While our result \eqref{mainbound} implies only H\"older continuity for $\|h\|_H\le1$, we can also obtain a locally Lipschitz continuity
\begin{equation}\label{mainbound:inhomversion}
\big\| {\widetilde{\sfM}}^{\infty;\eps}(\omega+h) : {\widetilde{\sfM}}^{\infty;\eps}(\omega) \big\|_{\textbf{\textsf{M}}(\scrV_{\varepsilon})_{w_a}} 
\leq \mathbb{L}(\omega) \max\big(\|h\|_H,\|h\|^{m({\bf B})}_H\big)
\end{equation}
for the inhomogeneous metric, by slightly modifying the above proof, which is consistent with Theorem 7 of \cite{Sch23}.

\medskip

\section{Large deviation principle for the BPHZ models}
\label{SectionLDP}

We make in this section the following two assumptions in addition to \textbf{Assumption A}.
\begin{itemize}
\item Every bounded set of $H$ is relatively compact in $\Omega$.
\item Write $\xi(\omega)=\omega$. As $\kappa\to0$, the family of random variables $\{\kappa\xi\}_{\kappa\in(0,1]}$ satisfies a large deviation principle with rate $\kappa^2$ and rate function $I:\Omega\to[0,\infty]$ given by
$$
I(\omega)=\begin{cases}
\frac12\|\omega\|_H^2&\text{if } \omega\in H,\\
\infty&\text{otherwise}.
\end{cases}
$$
That is one has
\begin{align*}
\limsup_{\kappa\to0}\kappa^2\log \bbP(\kappa\xi\in C)&\le -\inf_{\omega\in C} I(\omega)\qquad\forall\, C\subset \Omega:\text{closed},\\
\liminf_{\kappa\to0}\kappa^2\log \bbP(\kappa\xi\in O)&\ge -\inf_{\omega\in O} I(\omega)\qquad\forall\, O\subset \Omega:\text{open}.
\end{align*}
\end{itemize}
Under the same assumptions as in Theorem \ref{ThmConstructionBPHZ}, 
for every $\kappa\in(0,1]$, one can construct the BPHZ models ${\sf M}_\kappa^{n;\eps}$ on $\scrV_\eps$ lifted from $\kappa\xi*\varrho_n$, and the limit model $\widetilde{\sf M}_\kappa^{\infty;\eps}$. The aim of this section is to show the following theorem.

\ssk

\begin{thm} \label{thm:LDPforBPHZmodel}
As $\kappa\to0$, the $\textbf{\textsf{M}}(\scrV_{\varepsilon})_{w_a}$-valued random variables $\widetilde{\sf M}_\kappa^{\infty;\eps}$ satisfies a large deviation principle with rate $\kappa^2$ and rate function
$$
J({\sf M})=\inf\{I(h)\,;\,h\in H,\ L(h)={\sf M}\},
$$
where $L(h)\defeq{\sf M}^{h;\eps}$ is a naive model defined by setting ${\sf\Pi}^h(\ocircle)=h$. (Note that, unlike in the previous setting where $\ocircle$ and $\odot$ corresponded to $\xi$ and $h$, respectively, here $\ocircle$ corresponds to $h$.)
\end{thm}

\ssk

As above this means that one has
\begin{align*}
\limsup_{\kappa\to0}\kappa^2\log \bbP(\widetilde{\sf M}_\kappa^{\infty;\eps}\in C)&\le -\inf_{{\sf M}\in C} J({\sf M})\qquad\forall\, C\subset \textbf{\textsf{M}}(\scrV_{\varepsilon})_{w_a}:\text{closed},\\
\liminf_{\kappa\to0}\kappa^2\log \bbP(\widetilde{\sf M}_\kappa^{\infty;\eps}\in O)&\ge -\inf_{{\sf M}\in O} J({\sf M})\qquad\forall\, O\subset \textbf{\textsf{M}}(\scrV_{\varepsilon})_{w_a}:\text{open}.
\end{align*}   
Combining Theorem \ref{thm:LDPforBPHZmodel} with the continuity of the solution map $S:{\sf M}\mapsto u$ from a model $\sf M$ corresponding to the SPDE \eqref{eq:generalSPDE} to its maximal solution $u$ -- as stated in Theorem 7.8 and Corollary 7.12 of \cite{Hai14}, and more precisely in Theorems 2.21 and 5.21 of \cite{BCCH21} -- we obtain the following result.

\ssk

\begin{cor}
Assume that the equation is subcritical and that the assumptions of Theorem \ref{ThmConstructionBPHZ} are satisfied. For a fixed initial condition $u_0$, $u_{\kappa,n}=S(u_0,{\sf M}_\kappa^{n;\eps})$ is the maximal solution to the suitably renormalized equation of the form
$$
(\partial_t-\Delta)u_{\kappa,n} = \widetilde{\mcF}_{\kappa,n}(u_{\kappa,n},\nabla u_{\kappa,n};\kappa\xi*\varrho_n).
$$
As $n\to\infty$, the process $u_{\kappa,n}$ converges to $u_{\kappa,\infty}=S(u_0,\widetilde{\sf M}_\kappa^{\infty;\eps})$. As $\kappa\to0$, the family $u_{\kappa,\infty}$ satisfies a large deviation principle with rate $\kappa^2$ and rate function
\begin{align*}
K(u) &= \inf\{J({\sf M})\,;\,S(u_0,{\sf M})=u\} = \inf\{I(h)\,;\,S(u_0,L(h))=u\}   \\
        &= \inf\{I(h)\,;\,u(0,\cdot)=u_0,\ (\partial_t-\Delta)u = \mcF(u,\nabla u;h)\}.
\end{align*}
\end{cor}

\ssk

The key to the proof is Theorem \ref{Thm:Mn:Minfty} below, which is proved in the same way as Theorem \ref{ThmHolderEstimateBPHZ}. Since $\textbf{\textsf{M}}(\scrV_{\varepsilon})_{w_a}$ is not a linear space, however, we introduce an extended Polish space as follows.
This extended space is only used to state Theorem \ref{Thm:Mn:Minfty}.

\ssk

$\bullet$ Let $\sf M=(\Pi,g)$ be a pair consisting of a family of continuous linear maps $\big({\sf\Pi}_x:V\to C^{r_0,\mcQ}(w_a)\big)_{x\in\bbR^{1+d}}$ and a family of continuous linear functionals ${\sf g}=({\sf g}_x,{\sf g}_{yx})_{x,y\in\bbR^{1+d}}\subset(V_{\varepsilon}^+)^*$. 
(Although the object denoted by $\sf M=(\Pi,g)$ here is different from the one considered for models so far, we use the same notation.)
For such an $\sf M$, define
$$
\|{\sf M}\|_{\textbf{\textsf{M}}_{\textrm{hom}}(\scrV_{\varepsilon})_{w_a}} \defeq \max_{\tau\in{\bfB}\setminus\{X^k\}_k} \|{\sf\Pi}\res \tau\|_{\epsilon;w_a,\textrm{hom}}  +  \max_{\mu\in{\bf V}_{\varepsilon}^+\setminus\{X^k\}_k} \|{\sf g} \res \mu\|_{\epsilon;w_a,\textrm{hom}}
$$
in the same way as before.
We denote by $\widehat{\textbf{\textsf{M}}}_{\textrm{hom}}(\scrV_{\varepsilon})_{w_a}$ the completion, with respect to the above norm, of the set of all such $\sf M$ for which ${\sf\Pi}_x\tau(y)$, ${\sf g}_x(\mu)$, and ${\sf g}_{yx}(\mu)$ are smooth in $x,y$ for every $\tau,\mu$. Then $\widehat{\textbf{\textsf{M}}}_{\textrm{hom}}(\scrV_{\varepsilon})_{w_a}$ is a Polish space, and $\textbf{\textsf{M}}_{\textrm{hom}}(\scrV_{\varepsilon})_{w_a}$ is a closed subset of it.

\ssk

\begin{thm} \label{Thm:Mn:Minfty}
For all $\omega\in\Omega_0$ and $h\in H$, one has
 \begin{align*}
        \big\| \big({\widetilde{\sfM}}^{\infty;\eps}(\omega+h)-{\sfM}^{n;\eps}(\omega+h)\big) - \big( {\widetilde{\sfM}}^{\infty;\eps}(\omega) &- {\sfM}^{n;\eps}(\omega)\big)  \big\|_{\widehat{\textbf{\textsf{M}}}_{\textrm{\emph{hom}}}(\scrV_{\varepsilon})_{w_a}} \\
        &\leq \mathbb{L}_n(\omega) \max\big(\|h\|_H,\|h\|^{1/m({\bf B})}_H\big)
\end{align*}
with a random variable $\mathbb{L}_n(\omega)$ proportional to 
\[
\sup_{\|h\|_H\le1} \big\| {\overline{\sfM}}^{\infty,h;\eps}(\omega) \res {\sfM}^{n,h;\eps}(\omega) \big\|_{\textbf{\textsf{M}}_{\textrm{\emph{hom}}}(\scrU_{\varepsilon})_{w_a}},
\] 
which converges to $0$ in $L^p$ for any $p\ge1$.
\end{thm}

\ssk

\begin{cor} \label{Cor:ExpEquivl}
There exist a positive constant $C$ and a diverging sequence of positive constants $c_n$ such that 
\[
\bbE\big[\exp\big(c_n\|\widetilde{\sfM}^{\infty;\varepsilon}\res{\sf M}^{n;\eps}\|_{\textbf{\textsf{M}}_{\textrm{\emph{hom}}}(\scrV_{\varepsilon})_{w_a}}^2\big)\big] \le C.
\] 
\end{cor}

\ssk

\begin{Dem}
The proof is similar to that of Corollary \ref{cor:expintmodel}. Applying Proposition \ref{PropGeneralisedTCI} to $\textbf{\textsf{M}}=\widehat{\textbf{\textsf{M}}}_{\textrm{hom}}(\scrV_{\varepsilon})_{w_a}$, we see that the law $P_n$ of the random variable $\widetilde{\sf M}^{\infty;\eps}-{\sf M}^{n;\eps}$ satisfies a $(\ell_1(a_n(\cdot)),\|\cdot\|_{\widehat{\textbf{\textsf{M}}}_{\textrm{hom}}(\scrV_{\varepsilon})_{w_a}})$-transportation cost inequality, where $a_n$ is a positive constant proportional to $\|\mathbb{L}_n\|_{L^2(\bbP)}^{-1}$ diverging as $n\to\infty$.
Therefore, applying Theorem \ref{thm:expint} with $x_0={\bf0}$ and using that $\ell_1(t)\gtrsim t^2$ for $t\ge1$, we obtain
$$
\bbE\Big[\exp\Big(\tfrac{a_n^2}{10}
\big\{\big\|\widetilde{\sf M}^{\infty;\eps}-{\sf M}^{n;\eps}\big\|_{\widehat{\textbf{\textsf{M}}}_{\textrm{hom}}(\scrV_{\varepsilon})_{w_a}}
-\bbE\big[\big\|\widetilde{\sf M}^{\infty;\eps}-{\sf M}^{n;\eps}\big\|_{\widehat{\textbf{\textsf{M}}}_{\textrm{hom}}(\scrV_{\varepsilon})_{w_a}}\big]
\big\}^2\Big)
\Big]
\lesssim1.
$$
Since $b_n=\bbE\big[\big\|\widetilde{\sf M}^{\infty;\eps}-{\sf M}^{n;\eps}\big\|_{\widehat{\textbf{\textsf{M}}}_{\textrm{hom}}(\scrV_{\varepsilon})_{w_a}}\big]\to0$, by suitably redefining $a_n$ so that still $a_n\to\infty$ and $a_nb_n=O(1)$, we obtain the conclusion.
\end{Dem}

\ssk

Before proving Theorem \ref{thm:LDPforBPHZmodel}, note that ${\sf M}_\kappa^{n;\eps}$ admits the decomposition
$$
{\sf M}_\kappa^{n;\eps}=R_{\kappa,n}L(\kappa\xi*\varrho_n),
$$
where $R_{\kappa,n}$ is a map of the form \eqref{DefnRChi} with renormalization constants $\{\chi_{\kappa,n}(\tau)\}_{\tau\in{\bf B}_-}$.
Since $\chi_{\kappa,n}(\tau)$ is a homogeneous polynomial of $\kappa$ with degree $n_\circ(\tau)$, the number of symbol $\ocircle$ contained in $\tau$, it follows that
$$
\lim_{\kappa\to0}\chi_{\kappa,n}(\tau)=0.
$$
Therefore, $\lim_{\kappa\downarrow0}R_{\kappa,n}=\textrm{Id}$ for each $n\in\bbN$.
Define the continuous maps $F_{\kappa,n},F_{0,n}:\Omega\to\textbf{\textsf{M}}(\scrV_{\varepsilon})_{w_a}$ by
$$
F_{\kappa,n}(\omega)=R_{\kappa,n}L(\varrho_n*\omega),\qquad
F_{0,n}(\omega)=L(\varrho_n*\omega).
$$
Note that the following convergences hold.
\begin{itemize}
\item[(A)] For each $n\in\bbN$, $F_{\kappa,n}\to F_{0,n}$ uniformly over every bounded set of $\Omega$
\item[(B)] $F_{0,n}\to L$ uniformly over every bounded set of $H$.
\end{itemize}
(A) is obvious since $\varrho_n*\omega$ is bounded in $C^\infty$ on bounded subsets of $\Omega$. (B) is implied by the convergence ${\sf M}^{n,h_n;\eps}\to{\sf M}^{\infty,h;\eps}$ in $\textbf{\textsf{M}}(\scrU_{\varepsilon})_{w_a}$. Indeed, writing $F_{0,n}(h)=({\sf\Pi},{\sf g})$, one has
$$
{\sf\Pi}_x(\tau) = {\sf\Pi}_x^{n,h_n;\eps}(D_{\ocircle_\tau}\tau),\qquad
{\sf g}_{yx}(\sigma)={\sf g}_{yx}^{n,h_n;\eps}(D_{\ocircle_\sigma}\sigma).
$$

\ssk

The following argument is almost identical to that in Hairer \& Weber \cite{HW15}, but does not use Wiener chaos arguments. For the sake of completeness, we provide the proof.
We first recall the asymptotic contraction principle from Lemma 2.1.4 of \cite{DS89}.
We say that a function $I:S\to[0,\infty]$ on a metric space $S$ is a good rate function if it is lower semicontinuous, not identical to $\infty$, and $\{I\le\lambda\}$ is compact for any $\lambda\ge0$.

\ssk

\begin{lem}\label{lem:AsymptoticContractionPrinciple}
Let $(S,d)$ and $(T,\rho)$ be separable metric spaces and $\mathbb{I}:S\to[0,\infty]$ be a good rate function.
Let $\mathbb{F}:S\to T$ be a measurable map and assume that there exists a family $\{\mathbb{F}_\kappa:S\to T\}_{\kappa>0}$ of continuous maps such that $\mathbb{F}_\kappa\to \mathbb{F}$ uniformly over $\{\mathbb{I}\le\lambda\}$ for every $\lambda\in[0,\infty)$.
Then the maps
$$
\mathbb{J}_\kappa(t)=\inf\{\mathbb{I}(s)\,;\, s\in S,\ t=\mathbb{F}_\kappa(s)\},\qquad
\mathbb{J}(t)=\inf\{\mathbb{I}(s)\,;\, s\in S,\ t=\mathbb{F}(s)\}
$$
are good rate function on $T$.
Moreover, for any closed subset $C$ of $T$, one has
$$
\inf_{t\in C}\mathbb{J}(t)
=\lim_{\delta\to0}\liminf_{\kappa\to0}\inf_{t\in C_\delta}\mathbb{J}_\kappa(t),
$$
where $C_\delta=\{t\in T\,;\,\rho(t,C)\le\delta\}$.

In addition, assume that a family of Borel probability measures $\{\mu_\kappa\}_{\eps>0}$ on $S$ satisfies a large deviation principle with rate $\kappa^2$ and rate function $\mathbb{I}$, and
\begin{align}\label{lem:AsymptoticContractionPrinciple:Eq}
\limsup_{\kappa\to0}\kappa^2\log \mu_\kappa\big(\{s\in S\,;\,d(\mathbb{F}_\kappa(s),\mathbb{F}(s))>\lambda\}\big)=-\infty
\end{align}
for any $\lambda>0$.
Then the probability measures $\{\mu_\kappa\circ \mathbb{F}_\kappa^{-1}\}_{\kappa>0}$ on $T$ satisfy a large deviation principle with rate $\kappa^2$ and rate function $\mathbb{J}$.
\end{lem}

\ssk

\begin{cor}\label{cor:LDPn}
For each $n$, $\{{\sf M}_\kappa^{n;\eps}\}_{\kappa>0}$ satisfies a large deviation principle with rate $\kappa^2$ and rate function
$$
J_n({\sf M})=\inf\{I(h)\,;\,h\in H,\ F_{0,n}(h)={\sf M}\}
$$
\end{cor}

\ssk

\begin{Dem}
By the property (A), for any $\delta>0$, there exists an open neighborhood $O_\delta$ of $\{I\le\delta\}$ such that $F_{\kappa,n}\to F_{0,n}$ uniformly over $O_\delta$.
Therefore,
\begin{align*}
&\limsup_{\kappa\to0}\kappa^2\log\bbP(\| F_{\kappa,n}(\kappa\xi) \res F_{0,n}(\kappa\xi) \|_{\textbf{\textsf{M}}_{\textrm{hom}}(\scrV_{\varepsilon})_{w_a}}>\lambda)\\
&\le\limsup_{\kappa\to0}\kappa^2\log\bbP(\kappa\xi\in O_\delta^c)
\le-\inf_{\omega\in O_\delta^c}I(\omega)\le-\delta.
\end{align*}
Since $\delta>0$ is arbitrary, we obtain \eqref{lem:AsymptoticContractionPrinciple:Eq}.
\end{Dem}

\ssk

By the property (B) and by the first half part of Lemma \ref{lem:AsymptoticContractionPrinciple}, $J$ is a good rate function on $\textbf{\textsf{M}}(\scrV_{\varepsilon})_{w_a}$, and for any closed $C\subset\textbf{\textsf{M}}_{\text{\rm hom}}(\scrV_{\varepsilon})_{w_a}$, one has
\begin{equation}\label{lem:Inconvergence}
\inf_{{\sf M}\in C} J({\sf M})=\lim_{\delta\to0}\liminf_{n\to\infty}\inf_{{\sf M}\in C_\delta} J_n({\sf M}).
\end{equation}

\ssk

\begin{lem}\label{lem:ModelDilation}
Define the dilation operator $\delta_\kappa$ on $\textbf{\textsf{M}}(\scrV_{\varepsilon})_{w_a}$ by setting for $\sf M=(\Pi,g)$,
$$
\delta_\kappa({\sf\Pi})(\tau) = \kappa^{n_\circ(\tau)}{\sf\Pi}(\tau),\qquad
(\delta_\kappa{\sf g})(\mu)=\kappa^{n_\circ(\mu)}{\sf g}(\mu)
$$
for basis vectors.
Then one has $\|\delta_\kappa{\sf M}\res\delta_\kappa{\sf M}'\|_{\textbf{\textsf{M}}(\scrV_{\varepsilon})_{w_a}}=\kappa\|{\sf M}\res{\sf M}'\|_{\textbf{\textsf{M}}(\scrV_{\varepsilon})_{w_a}}$ and
$$
\delta_\kappa{\sf M}_1^{n;\eps}={\sf M}_\kappa^{n;\eps},\qquad
\delta_\kappa\widetilde{\sf M}_1^{\infty;\eps}=\widetilde{\sf M}_\kappa^{\infty;\eps}.
$$
\end{lem}

\ssk

\begin{Dem}
In the construction of the naive model and the BPHZ model, ${\sf\Pi}(\tau)$ and ${\sf g}(\mu)$ are always kept as homogeneous polynomials in the noise $\xi$, so the construction ensures the result.
\end{Dem}

\ssk

\begin{Dem}[{\bf of Theorem \ref{thm:LDPforBPHZmodel}}]
{\bf Upper bound}: Let $C\subset \textbf{\textsf{M}}_{\textrm{hom}}(\scrV_{\varepsilon})_{w_a}$ be an arbitrary closed set. For any $\lambda>0$, we decompose
$$
\bbP(\widetilde{\sf M}_\kappa^{\infty;\varepsilon}\in C)
\le \bbP({\sf M}_\kappa^{n;\eps}\in C_\lambda)+\bbP\big(\big\| {\widetilde{\sfM}}_\kappa^{\infty;\eps} \res {\sfM}_\kappa^{n;\eps} \big\|_{\textbf{\textsf{M}}_{\textrm{hom}}(\scrV_{\varepsilon})_{w_a}}>\lambda\big).
$$
By Lemma \ref{lem:ModelDilation} and Corollary \ref{Cor:ExpEquivl},
\begin{align*}
\bbP\big(\big\| {\widetilde{\sfM}}_\kappa^{\infty;\eps} \res {\sfM}_\kappa^{n;\eps} \big\|_{\textbf{\textsf{M}}_{\textrm{hom}}(\scrV_{\varepsilon})_{w_a}}>\lambda\big)
&=\bbP\big(\big\| {\widetilde{\sfM}}_1^{\infty;\eps} \res {\sfM}_1^{n;\eps} \big\|_{\textbf{\textsf{M}}_{\textrm{hom}}(\scrV_{\varepsilon})_{w_a}}>\lambda/\kappa\big)\\
&\lesssim e^{-c_n\lambda^2/\kappa^2}.
\end{align*}
Thus by Corollary \ref{cor:LDPn}, we have
\begin{align*}
\limsup_{\kappa\to0}\kappa^2\log\bbP(\widetilde{\sf M}_\kappa^{\infty;\varepsilon}\in C)
&\le \bigg(\limsup_{\kappa\to0}\kappa^2\log\bbP({\sf M}_\kappa^{n;\varepsilon}\in C_\lambda)\bigg)\vee(-c_n\lambda^2)\\
&\le \bigg(-\inf_{{\sf M}\in C_\lambda} J_n({\sf M})\bigg)\vee(-c_n\lambda^2).
\end{align*}
Taking the limit $n\to\infty$, we have
$$
\limsup_{\kappa\to0}\kappa^2\log\bbP(\widetilde{\sf M}_\kappa^{\infty;\varepsilon}\in C)
\le -\liminf_{n\to\infty}\inf_{{\sf M}\in C_\lambda} J_n({\sf M}).
$$
Taking the limit $\lambda\to0$, we have the desired upper bound by \eqref{lem:Inconvergence}.

\medskip

\noindent
{\bf Lower bound}: We fix a nonempty open set $O\subset \textbf{\textsf{M}}_{\textrm{hom}}(\scrV_{\varepsilon})_{w_a}$, an element ${\sf M}\in O$ and $\lambda>0$ such that $B_{\textrm{hom}}({\sf M},2\lambda)\subset O$, the closed ball centered at $\sf M$ with radius $2\lambda$ with respect to the metric $\|\cdot,\cdot\|_{\textbf{\textsf{M}}_{\textrm{hom}}(\scrV_{\varepsilon})_{w_a}}$. By decomposing
$$
\bbP\big({\sf M}_\kappa^{n;\eps}\in B_{\textrm{hom}}({\sf M},\lambda)\big)
\le \bbP(\widetilde{\sf M}_\kappa^{\infty;\eps}\in O)+\bbP\big(\big\| {\widetilde{\sfM}}_\kappa^{\infty;\eps} \res {\sfM}_\kappa^{n;\eps} \big\|_{\textbf{\textsf{M}}_{\textrm{hom}}(\scrV_{\varepsilon})_{w_a}}>\lambda\big),
$$
we have
\begin{align*}
-\inf_{{\sf M}'\in B_{\textrm{hom}}({\sf M},\lambda/2)} J_n({\sf M}')
&\le\liminf_{\kappa\to0}\kappa^2\log \bbP\big({\sf M}_\kappa^{n;\eps}\in B_{\textrm{hom}}({\sf M},\lambda)\big)\\
&\le\Big(\liminf_{\kappa\to0}\kappa^2\log \bbP(\widetilde{\sf M}_\kappa^{\infty;\eps}\in O)\Big)\vee
(-c_n\lambda^2).
\end{align*}
By first taking the limit $n\to\infty$ and then letting $\lambda\to0$ as before, we have
\begin{align*}
\liminf_{\kappa\downarrow0}\kappa^2\log \bbP(\widetilde{\sf M}_\kappa^{\infty;\eps}\in O)
\ge -J({\sf M}).
\end{align*}
by \eqref{lem:Inconvergence}. Since ${\sf M}\in O$ is arbitrary, we obtain the desired lower bound.
\end{Dem}

\medskip

\appendix   

\section{A relation between transportation cost and spectral gap inequalities}
\label{SectionAppendix}

Denote by $\iota : H\rightarrow \Omega$ the continuous injection of $H$ into $\Omega$.

\ssk

\begin{lem} \label{lem:expintegrability}
Under Assumption A one has for all small enough $\lambda>0$
$$
\int_\Omega e^{\lambda\|\omega\|^2}\bbP(d\omega) < \infty.
$$
\end{lem}

\ssk

\begin{Dem}
Following the arguments of the proofs of Proposition 1.6 and Theorem 1.7 in Gozlan \& L\'eonard's review \cite{GL10}, we have the existence of some positive finite constants $a_1,a_2$ such that one has 
\begin{equation}\label{eq:lem:expintegrability}
\bbP(F-m_F>r)\le a_1e^{-a_2 r^2}
\end{equation}
for all $r>0$ and any $1$-Lipschitz function $F:\Omega\to\bbR$ with respect to the extended distance $d_H$, with $m_F$ is the median of $F$. One gets the conclusion by applying the inequality \eqref{eq:lem:expintegrability} to the function $F(\omega)=\|\omega\|_\Omega/\|\iota\|_{H\to\Omega}$.
\end{Dem}

\ssk

\begin{thm}
A probability measure satisfying Assumption A also satisfies the spectral gap inequality \eqref{EqSG} with $a_1=1/(2a_\circ)$.
\end{thm}

\ssk

\begin{Dem}
Let $\mu$ be a Borel probability measure on $\Omega$ satisfying Assumption A, and let $F$ be a cylindrical function of the form $F(\omega) = f\big(\varphi_1(\omega),\dots,\varphi_n(\omega)\big)$ with $\varphi_1,\dots,\varphi_n\in\Omega^*$ and at most polynomially growing $f\in C^\infty(\bbR^n)$. By applying the Gram--Schmidt process, we may assume that $\{\iota^*(\varphi_k)\}_{k=1}^n\subset H^*$ is an orthonormal system in $H^*$.
Then note that the projection map $p_n : \Omega\to\bbR^n$ defined by $p_n(\omega) = (\varphi_1(\omega),\dots,\varphi_n(\omega))$ is a contraction with respect to the extended distance $d_H=\sqrt{c_H}$ on $\Omega$ and the Euclidean distance $d_n$ on $\bbR^n$ as
$$
d_n\big(p_n(\omega+\iota(h)) , p_n(\omega)\big)^2 = \sum_{k=1}^n \big| \varphi_k(\iota(h)) \big|^2 = \sum_{k=1}^n \big| \iota^*(\varphi_k)(h) \big|^2 \leq \|h\|_{H}^2=d_H(\omega+\iota(h),\omega)^2.
$$
Let $\mu_n=\mu\circ p_n^{-1}$ be the image measure of $\mu$ on $\bbR^n$ by the application $p_n$. For any probability measure $\nu_n(dx)=g_n(x)\mu_n(dx)$ on $\bbR^n$ with $g_n\ge0$, define a probability measure $\rho_n(d\omega) \defeq g_n(p_n(\omega))\mu(d\omega)$ on $\Omega$. Then since $\frac{d\rho_n}{d\mu}(\omega)=\frac{d\nu_n}{d\mu_n}(p_n(\omega))$, we have $\mcH(\nu_n|\mu_n)=\mcH(\rho_n|\mu)$. Moreover, for any coupling $\Theta_n$ of $\rho_n$ and $\mu$, by setting $\theta_n \defeq \Theta_n \circ (p_n\times p_n)^{-1}$ we have
\begin{align*}
\mcW_{d_n^2}(\nu_n,\mu_n) \leq \int_{\bbR^n\times\bbR^n} d_n(x,y)^2 \theta_n(dxdy) &= \int_{\Omega\times\Omega} d_n(p_n(\omega_1),p_n(\omega_2))^2 \Theta_n(d\omega_1d\omega_2)   \\
&\leq \int_{\Omega\times\Omega}d_H(\omega_1,\omega_2)^2 \, \Theta_n(d\omega_1d\omega_2),
\end{align*}
which implies $\mcT_{d_n^2}(\nu_n,\mu_n) \leq \mcT_{d_H^2}(\rho_n,\mu)$. Hence $\mu_n$ satisfies a $(a_\circ \textrm{Id},d_n^2)$-transportation cost inequality. Then by Otto and Villani's theorem \cite{OV00} -- see also \cite[Proposition 8.11]{GL10}), the probability $\mu_n$ satisfies the spectral gap inequality
$$
2a_\circ \bbE_{\mu_n}[|f-\bbE_{\mu_n}[f]|^2]\le \bbE_{\mu_n}[|\nabla f|^2],
$$
where $\nabla f=(\partial_1f,\cdots,\partial_nf)$ is the gradient vector. We note that we have $\bbE_{\mu_n}[|f-\bbE_{\mu_n}[f]|^2]=\bbE_{\mu}[|F-\bbE_{\mu}[F]|^2]$ by definition of $\mu_n$. By taking a dual system $\{e_k\}_{k=1}^n\subset H$ of $\{\iota^*(\varphi_k)\}\subset H^*$ such that $\iota^*(\varphi_k)(e_\ell)={\bf1}_{k=\ell}$, we have
$$
|\nabla f|^2(p_n(\omega)) = \sum_{k=1}^n|\partial_kf(p_n(\omega))|^2 = \sum_{k=1}^n |d_\omega F(e_k)|^2 \le\|d_\omega F\|_{H^*}^2.
$$
This shows the spectral gap inequality $2a_\circ \bbE_{\mu}[|F-\bbE_{\mu}[F]|^2]\le \bbE_{\mu}[\|d_\omega F\|_{H^*}^2]$.

\medskip

We actually use the spectral gap inequality for functionals of the form $\mcQ_t(x,{\sf\Pi}_x^{n;\eps}\tau)$. These functionals are not cylindrical functions but they have the form
\begin{equation} \label{eq:extendedcylindrical}
F(\omega) = \int_{\bbR^n}f(\varphi_x(\omega))\theta(dx)
\end{equation}
for a continuous bounded function $x\mapsto\varphi_x$ taking values in the space of bounded linear operators $\Omega\to\bbR^n$, for an at most polynomially growing function $f\in C^\infty(\bbR^n)$ and a finite signed Borel measure $\theta$ on $\bbR^n$. Functions ${\sf\Pi}_x^{n;\eps}\tau$ are generated from regularized functions $\varrho*\omega$ with some smooth mollifier $\varrho$ and by applying the three kinds of operations inductively; (1) linear combinations, (2) multiplications, and (3) integral operator $\mcK$. Therefore, functionals $\mcQ_t(x,{\sf\Pi}_x^{n;\eps}\tau)$ can be written in the form \eqref{eq:extendedcylindrical}.

\ssk

Let $F$ be of the form \eqref{eq:extendedcylindrical}. By the bounded continuity of the map $x\mapsto f(\varphi_x(\omega))$, there is a finite partition $\{A_{m,k}\}_{k=1}^{N_m}$ of $\bbR^n$ and a set of reference points $\{a_{m,k}\in A_{m,k}\}_{k=1}^{N_m}$ for each $m\in\bbN$ such that
$$
F_m(\omega)\defeq\sum_{k=1}^{N_m}f(\varphi_{a_{m,k}}(\omega))\theta(A_{m,k})
$$
converges to $F(\omega)$ as $m$ goes to $\infty$ for all $\omega\in\Omega$. Moreover, $d_\omega F_m(h)$ also converges to $d_\omega F(h)$ uniformly over $\|h\|_H\le1$ for all $\omega\in\Omega$. Since $F_m$ is a cylindrical function as discussed before, it satisfies the inequality
$$
2a_\circ \bbE_{\mu}[|F_m-\bbE_{\mu}[F_m]|^2]\le \bbE_{\mu}[\|d_\omega F_m\|_{H^*}^2] = \bbE\bigg[\sup_{\Vert h\Vert_H\leq 1}|d_\omega F_m(h)|^2\bigg].
$$
By the pointwise convergence and the uniform integrabilities of $F_m$ and $d_\omega F_m$ given by Lemma \ref{lem:expintegrability}, we can take the limit $m\to\infty$ in the above inequality and obtain the conclusion.
\end{Dem}

\vspace{0.5cm}

\noindent \textcolor{gray}{$\bullet$} {\sf Bailleul I.} -- Univ Brest, CNRS, LMBA - UMR 6205, F-29238 Brest, France.   \\
{\it E-mail}: ismael.bailleul@univ-brest.fr

\medskip

\noindent \textcolor{gray}{$\bullet$} {\sf Hoshino M.} --  Department of Mathematics, Institute of Science Tokyo, Japan.   \\
{\it E-mail}: hoshino@math.sci.isct.ac.jp

\medskip

\noindent \textcolor{gray}{$\bullet$} {\sf Takano R.} --  The University of Osaka, Japan.   \\
{\it E-mail}: rtakano@sigmath.es.osaka-u.ac.jp 


\bigskip


\begin{thebibliography}{99}



\bibitem{BH23}
{\sc Bailleul I.} and {\sc Hoshino M.},
\newblock {\em Random models on regularity-integrability structures}.
\newblock arXiv:2310.10202.

\bibitem{BHRSGuide}
{\sc Bailleul I.} and {\sc Hoshino M.},
\newblock {\em A tourist's guide to to regularity structures and singular stochastic PDEs}.
\newblock doi:10.4171/EMSS/87. To appear in EMS Surv. Math. Sci. (2025+).


\bibitem{BauerschmidtBodineau}
{\sc Bauerschmidt R.} and {\sc Bodineau T.},
\newblock {\em Log-Sobolev inequality for the continuum sine-Gordon model}.
\newblock Comm. Pure Appl. Math. {\bf 74}(10):2064--2113, (2021).

\bibitem{BauerschmidtDagallier}
{\sc Bauerschmidt R.} and {\sc Dagallier B.},
\newblock {\em Log-Sobolev inequality for the $\Phi^4_2$ and $\Phi^4_3$}.
\newblock  Comm. Pure Appl. Math. {\bf 77}(5):2579--2612, (2024).

\bibitem{Bru18}
{\sc Bruned Y.}, 
\newblock {\em Recursive formulae in regularity structures}.
\newblock Stoch PDE: Anal Comp {\bf 6}:525--564 (2018).

{\color{red}
\bibitem{BCCH21}
{\sc Bruned Y., Chandra A., Chevyrev I.}, and {\sc Hairer M.},
\newblock {\em Renormalising SPDEs in regularity structures}.
\newblock J. Eur. Math. Soc. (JEMS) {\bf 23} (2021), 869--947.
}

\bibitem{BHZ}
{\sc Bruned Y., Hairer M. {\rm and} Zambotti L.},
\newblock {\em Algebraic renormalisation of regularity structures}.
\newblock Invent. Math. {\bf 215}(3):1039--1156 (2019).


\bibitem{ChandraHairer}
{\sc Chandra A.} and {\sc Hairer M.},
\newblock {\em An analytic BPHZ theorem for Regularity Structures}.
\newblock arXiv:1612.08138 (2016).

\bibitem{CMW23}
{\sc Chandra A., Moinat A. {\rm and} Weber H.},
\newblock {\em A Priori Bounds for the  Equation in the Full Sub-critical Regime}. 
\newblock Arch Rational Mech Anal {\bf 247}, 48 (2023).

\bibitem{CFW24}
{\sc Chandra A., de Lima Feltes G. {\rm and} Weber H.}, 
\newblock {\em A priori bounds for 2-d generalised Parabolic Anderson Model}.
\newblock arXiv:2402.05544 (2024).

\bibitem{ChevyrevGubinelli}
{\sc Chevyrev I. {\rm and} Gubinelli M.},
\newblock {\em Large field problem in coercive singular PDEs}.
\newblock arXiv:2510.20716 (2025).

\bibitem{DL22}
{\sc Dai Y.} and {\sc Li R.N.}, 
\newblock {\em Transportation Inequalities for Stochastic Heat Equation with Rough Dependence in Space}.
\newblock Acta. Math. Sin.-English Ser. {\bf 38}:2019--2038 (2022).

\bibitem{DS89}
{\sc Deuschel J.-D.} and {\sc Stroock D. W.},
\newblock {\em Large deviations}.
\newblock Pure Appl. Math., {\bf 137}, Academic Press, Inc., Boston, MA (1989).

\bibitem{DGW04}
{\sc Djellout H., Guillin A. {\rm and} Wu L.}, 
\newblock {\em Transportation cost-information inequalities and applications to random dynamical systems and diffusions}.
\newblock Ann. Probab. {\bf 32}:2702--2732 (2004).

\bibitem{EsquivelWeber}
{\sc Esquivel S.} and {\sc Weber H.},
\newblock {\em A priori bounds for the dynamic fractional $\Phi^4$ model on $\mathbb{T}^3$ in the full subcritical regime}.
\newblock arXiv:2411.16536 (2024).

\bibitem{FU02}
{\sc Feyel D. {\rm and} \"Ust\"unel A. S.}, 
\newblock {\em Measure transport on Wiener space and the Girsanov theorem}. 
\newblock Comptes Rendus Mathematique {\bf 334}:1025--1028 (2002).


\bibitem{GJ23}
{\sc Gasteratos I. {\rm and} Jacquier A.},
\newblock {\em Transportation-cost inequalities for non-linear Gaussian functionals}.
\newblock arXiv:2310.05750 (2023).

\bibitem{Goz06}
{\sc Gozlan N.},
\newblock {\em Integral criteria for transportation cost inequalities}.
\newblock Elect. Comm. in Probab. {\bf 11}:64--77 (2006).


\bibitem{GL10}
{\sc Gozlan N. {\rm and} L\'{e}onard C.}, 
\newblock {\em Transport inequalities. A survey}.
\newblock Markov Proc. Rel. Fields {\bf 16}:635--736 (2010).

\bibitem{GH21}
{\sc Gubinelli M. {\rm and} Hofmanov\'a M.},
\newblock {\em A PDE Construction of the Euclidean $\Phi^4_3$ quantum field theory}.
\newblock Comm. Math. Phys. {\bf 384}:1--75 (2021).


\bibitem{Hai14}
{\sc Hairer M.}, 
\newblock {\em A theory of regularity structures},
\newblock Invent. Math. {\bf 198}:269--504 (2014).

\bibitem{Hai24}
{\sc Hairer M.},
\newblock {\em Renormalisation in the presence of variance blowup}.
\newblock Ann. Probab. {\bf 53}(5):1958--1985 (2025).

\bibitem{HairerMattingly}
{\sc Hairer M.} and {\sc Mattingly, J.},
\newblock {\em The strong Feller property for singular stochastic PDEs}.
\newblock Ann. Inst. H. Poincar\'e {\bf 54}(3):1314--1340 (2018).


\bibitem{HairerSchonbauer}
{\sc Hairer M.} and {\sc Sch\"onbauer P.},
\newblock {\em The support of singular stochastic partial differential equations}.
\newblock Forum Math. Pi {\bf 10}(e1):1--127 (2022).

\bibitem{HS22}
{\sc Hairer M. {\rm and} Steele R.},
\newblock {\em The $\Phi^4_3$ measure has sub-Gaussian tails}. 
\newblock J. Stat. Phys. {\bf 186}(3), Paper No. 38, 25 pp (2022).

\bibitem{HS23}
{\sc Hairer M. {\rm and} Steele R.},
\newblock {\em The} BPHZ {\em Theorem for Regularity Structures via the Spectral Gap Inequality}.
Arch. Rat. Mech. Anal. {\bf 248}, 9 (2024).

\bibitem{HW15}
{\sc Hairer M. {\rm and} Weber H.},
\newblock {\em Large deviations for white-noise driven, nonlinear stochastic PDEs in two dimensions}.
\newblock Ann. Fac. Sci. Toulouse Math. (6) {\bf 24}:55--92, (2015).

\bibitem{SemigroupMasato}
{\sc Hoshino M.},
\newblock {\em A semigroup approach to reconstruction theorems and multilevel Schauder estimates}.
\newblock Annales Henri Lebesgue, {\bf 8}:151--180, (2025).


\bibitem{KS19}
{\sc Khoshnevisan D. {\rm and} Sarantsev A.}, 
\newblock {\em Talagrand concentration inequalities for stochastic partial differential equations}.
\newblock Stoch. PDE: Anal. Comp. {\bf 7}:679--698 (2019).


\bibitem{LW24}
{\sc Li R. {\rm and} Wang X.},
\newblock {\em Talagrand's transportation inequality for SPDEs with locally monotone drifts}.
\newblock Statistics \& Probability Letters {\bf 204}:109945 (2024).

\bibitem{LOTT}
{\sc Linares P., Otto F., Tempelmayr M. {\rm and} Tsatsoulis P.},
\newblock {\em A diagram-free approach to the stochastic estimates in regularity structures}.
\newblock Invent. math. {\bf 237}:1469--1565, (2024).

\bibitem{Lyons98}
{\sc Lyons T.J.},
\newblock {\em Differential equations driven by rough signals}.
\newblock Rev. Mat. Iberoamericana {\bf 14}(2):215--310 (1998).

\bibitem{Mar86}
{\sc Marton K.},
\newblock {\em A simple proof of the blowing-up lemma}. 
\newblock IEEE Trans. Inform. Theory, {\bf 32}(3):445--446 (1986).


\bibitem{MW20}
{\sc Moinat A. {\rm and} Weber H.},
\newblock {\em Space-time localisation for the dynamic $\Phi^4_3$ model}.
\newblock Comm. Pure Appl. Math. {\bf 73}:2519--2555 (2020).

\bibitem{OV00}
{\sc Otto F. and Villani C.}, 
\newblock {\em Generalization of an inequality by Talagrand and links with the logarithmic Sobolev inequality}.
\newblock J. Funct. Anal. {\bf 173}:361--400 (2000).

\bibitem{Rie17}
{\sc Riedel S.}, 
\newblock {\em Transportation-cost inequalities for diffusions driven by Gaussian processes}.
\newblock Electronic Journal of Probability, {\bf 22} (2017), pp. 1--26.

\bibitem{Sau12}
{\sc Saussereau B.}, 
\newblock {\em Transportation inequalities for stochastic differential equations driven by a fractional Brownian motion}. 
\newblock Bernoulli, (2012), pp. 1--23.

\bibitem{Sch23}
{\sc Sch\"onbauer, P.},
\newblock {\em Malliavin calculus and densities for singular stochastic partial differential equations}.
\newblock Probab. Theory Related Fields {\bf 186} (2023), 643--713.

\bibitem{SW20}
{\sc Shang S. {\rm and} Wang R.}, 
\newblock {\em Transportation Inequalities Under Uniform Metric for a Stochastic Heat Equation Driven by Time-White and Space-Colored Noise}.
\newblock Acta Appl Math {\bf 170}, 81-97 (2020).

\bibitem{SZ19}
{\sc Shang S. {\rm and} Zhang T.}, 
\newblock {\em Talagrand concentration inequalities for stochastic heat-type equationsunder uniform distance}.
\newblock Electron. J. Probab., {\bf 24}:(129), 1-15, (2019). 

\bibitem{SZZ}
{\sc Shen H., Zhu, R.} and {\sc Zhu Z.},
\newblock {\em Global well-posedness for $2d$ generalized parabolic Anderson model via paracontrolled calculus}.
\newblock Stoch. PDEs: Anal. and Comput. (2025)

\bibitem{Tal96}
{\sc Talagrand M.}, 
\newblock {\em Transportation cost for Gaussian and other product measures}.
\newblock Geometric \& Functional Analysis {\bf 6} (1996), pp. 587--600.


\bibitem{WZ06}
{\sc Wu L. {\rm and} Zhang Z.},
\newblock {\em Talagrand's $T_2$-transportation inequality and log-Sobolev inequality for dissipative SPDEs and applications to reaction-diffusion equations}. 
\newblock Chinese Ann. Math. Ser. B, {\bf 27}:(3), 243-262, 2006.
\end{thebibliography}
\end{document}